\documentclass[12pt,a4,paper]{article}
\usepackage{amssymb,bm}
\usepackage{tipa,undertilde}
\renewcommand{\thefootnote}

\newtheorem{sect}{\S}
\newtheorem{defn}{Definition}[sect]
\newtheorem{cla}[defn]{Claim}
\newtheorem{prop}[defn]{Proposition}
\newtheorem{theo}[defn]{Theorem}

\begin{document}

\begin{center}
{\Large \bf{Mad Spectra}}\\
\vspace{.8cm} Saharon Shelah and Otmar Spinas\par \vspace{.8cm}
\end{center}

\begin{abstract}
\noindent The mad spectrum is the set of all cardinalities of infinite maximal
almost disjoint families on $\omega$. We treat the problem to
characterize those sets $\mathcal{A}$ which, in some forcing extension
of the universe, can be the mad spectrum. We solve this problem to some
extent. What remains open is the possible values of $\min (\mathcal{A})$
and $\max (\mathcal{A})$.
\end{abstract}\vspace{-1.5cm}

\footnote{1991 Mathematics Subject Classification. Primary: 03E35, 03E17; Key words
and phases: set theory, forcing, MAD families, set theory of the reals.\vspace{1ex}}

\footnote{The first author would like to thank the Israel Science Foundation
(ISF grant no. 1053/11) and the National Science Foundation (NSF grant no. DMS 1101597)
for partial support of this research. Publication 1038.\vspace{1ex}}

\footnote{The second author thanks the Deutsche Forschungsgemeinschaft (DFG) for
partial support (grant no. SP 683/3-1).}

{\bf \S\, 0 Introduction}\\

Recall that $A \subseteq [\omega]^\omega$ is called \textbf{almost disjoint}
(\textbf{a. d.} for short), if $a \cap b$ is finite for all $a, b \in A, a \neq b$.
Such $A$ is called \textbf{maximal almost disjoint} (\textbf{mad} for short), if
it is maximal with respect to $\subseteq$ among a. d. families. An easy diagonalization
shows that every infinite mad family is uncountable. The well-known cardinal invariant
$\frak{a}$ is defined as the minimal cardinality of an infinite mad family. Over the
past decades, much work has been done to understand this cardinal. We only mention
[Sh700] and [Br]. In the first one, the consistency of $\frak{d} < \frak{a}$ with
ZFC was proved, in the second one, which further develops the ideas of [Sh700], it
was shown that consistently $a = \aleph_\omega$.

It is natural to study mad families in more general ways, e. g. investigate the
\textbf{mad spectrum}, i. e. the set $\mathcal{A}$ of all infinite cardinals that
are the cardinality of some mad family. This problem has been attacked already
in the early period of forcing by Hechler [H]. There are two obvious restrictions
$\mathcal{A}$ must satisfy. Firstly, $\mathcal{A}$ has $2^{\aleph_0}$ as its
maximum, and, secondly, $\mathcal{A}$ is closed under singular limits (see [H,
Theorem 3.1]). For the first one, notice that there are always a. d. and hence
mad families of size $2^{\aleph_0}$. For the second one, if $\mu$ is a singular
limit of $\mathcal{A}$, say $\mu = \mathop{\Sigma}\limits_{i < cf (\mu)} \mu_i$ with
$cf (\mu) \leqslant \mu_i < \mu$ and $\mu_i \in \mathcal{A}$ for all $i < cf (\mu)$,
choose mad families $B_i$ with $|B_i| = \mu_i$. Let $B_0 = \{b_\nu : \nu < \mu_0 \}$
and fix bijections $\pi_\nu : \omega \rightarrow b_\nu, \nu < \mu_0$. For $i < cf (\mu)$
let $B'_{1+i} = \{\pi_i [b] : b \in B_{1+i}\}$. Then $(B_0 \setminus \{b_\nu : \nu <
cf (\mu)\}) \cup \bigcup\limits_{i < cf (\mu)} B'_{1+i}$ is a mad family of size
$\mu$.

It is natural to try to characterize those sets $\mathbb{C} \subseteq \mathbf{Card}$
which are the mad spectrum of some forcing extension of $\mathbf{V}$. Under the
assumption $\mathbf{V} \models GCH$ Hechler has constructed some c. c. c. forcing
notion $\mathbb{P}$, such that $\Vdash_\mathbb{P} \mathbb{C} = \utilde{\mathcal{A}}$
($\utilde{\mathcal{A}}$ is a $\mathbb{P}$-name for the mad spectrum), provided that:

\begin{itemize}
\item[(a)] $\mathbb{C}$ is a set of uncountable cardinals;
\item[(b)] $\mathbb{C}$ is closed under singular limits;
\item[(c)] if $\mu \in \mathbb{C}$ has $cf (\mu) = \aleph_0$, then $\mu = \sup
(\mathbb{C} \cap \mu)$;
\item[(d)] $\max (\mathbb{C})$ exists and $\max (\mathbb{C})^{\aleph_0} = \max
(\mathbb{C})$;
\item[(e)] $\aleph_1 \in \mathbb{C}$;
\item[(f)] if $\mu \in \mathbb{C}{\rm ard}$ and $\aleph_1 < \mu \leqslant |\mathbb{C}|$, then
$\mu \in \mathbb{C}$;
\item[(g)] if $\mu \in \mathbb{C}, cf (\mu) = \aleph_0$, then $\mu^+ \in
\mathbb{C}$.
\end{itemize}

The question remained open whether (c), (e), (f), (g) are necessary assumptions.
In particular, Raghavan has asked whether consistently $\aleph_\omega \in  \mathcal{A}$
but $\aleph_{\omega + 1} \notin \mathcal{A}$.

In this paper we show that for every  $\mathbb{C} \subseteq \mathbf{Card}$ with
properties (a), (b), (c), (d) there exists some c. c. c. forcing $\mathbb{P}_\mathbb{C}$
with $\Vdash_{\mathbb{P}_\mathbb{C}} \mathbb{C} = \utilde{\mathcal{A}}$, provided that
$\vartheta := \min (\mathbb{C})$ satisfies $\vartheta = \vartheta^{< \vartheta}$ (hence
$\vartheta$ regular) and $\max (\mathbb{C})^{< \vartheta} = \max (\mathbb{C})$. In particular,
we answer Raghavan's question positively. By Brendle's result mentioned above this is not
a complete characterization of the possible mad spectra.

{\sect {\bf The obvious Forcing}}

\begin{defn}\label{coco} We say $\mathbb{C}$ is a \textbf{potential mad spectrum}
(p.m.s. for short), if the following hold:

\begin{itemize}
\item[(a)] $\mathbb{C}$ is a set of cardinals;
\item[(b)] $\min (\mathbb{C}) \geqslant \aleph_1$;
\item[(c)] $\max (\mathbb{C})$ exists;
\item[(d)] $\max (\mathbb{C})^{\aleph_0} = \max (\mathbb{C})$;
\item[(e)] if $\mu$ is singular and $\mu = \sup (\mathbb{C} \cap \mu)$,
then $\mu \in \mathbb{C}$;
\item[(f)] as a technical assumption we ask that $\mathbb{C}$ has $\max
(\mathbb{C})$ as a member $\max (\mathbb{C})$ times, and we write them
as $\langle \Upsilon_i : i < \max (\mathbb{C})\rangle$.
\end{itemize}
\end{defn}

\begin{defn}\label{hei}
For any uncountable cardinal $\mu$ let $\mathbb{Q}_\mu$ be the following
forcing notion:

\begin{itemize}
\item[(A)] $p \in \mathbb{Q}_\mu$ iff for some unique $u = {\rm dom} (p)$ and $n = n_p < \omega$
we have
\begin{itemize}
\item[(a)]  $u \subseteq \mu$ is finite;
\item[(b)]  $p : u \rightarrow {^n2}$;
\end{itemize}

\item[(B)] $\mathbb{Q}_\mu \models p \leqslant q$ iff

\begin{itemize}
\item[(a)]  $p, q \in \mathbb{Q}_\mu$;
\item[(b)]  ${\rm dom} (p) \subseteq {\rm dom} (q)$;
\item[(c)]  if $\alpha \in {\rm dom} (p)$, then $p (\alpha) \unlhd q (\alpha)$ (hence $n_p \leqslant n_q$);
\item[(d)]  if $\alpha, \beta \in {\rm dom} (p), \alpha \neq \beta$ and $n \in [n_p, n_q)$, then
$q (\alpha) (n) = 0$ or $q (\beta) (n) = 0$.
\end{itemize}
\end{itemize}
\end{defn}

\textbf{Remark}
Note that if $p, q \in \mathbb{Q}_\mu$ with $n_p \leqslant n_q$ are incompatible,
then there exist $\alpha \in u_p \cap u_q$ and $n < n_p$, so that $p (\alpha) (n)
\neq q (\alpha) (n)$, or else $n_p < n_q$ and there exist $\alpha, \beta \in
u_p \cap u_q$ and $n \in [n_p, n_q)$ so that $\alpha \neq \beta$ and $q
(\alpha) (n) = q (\beta) (n) = 1$.

Recall that a forcing has the \textbf{Knaster} property, if every uncountable subset
has an uncountable subset such that any two of its elements are compatible.
Applying the $\Delta$-system lemma we easily get:

\begin{cla}
Forcing $\mathbb{Q}_\mu$ has the Knaster property, hence is c.c.c.
\end{cla}

\begin{defn}
1) For $u \subseteq \mu$ let $\mathbb{Q}_{\mu, u}$ be the forcing $\mathbb{Q}_\mu$ restricted
to $\{p \in \mathbb{Q}_\mu : {\rm dom} (p) \subseteq u\}$.\\
2) If $\mathbf{x} = (\mu_1, u_1, \mu_2, u_2, h)$ is such that $u_l \subseteq \mu_l$ for
$l = 1,2$ and $h$ is a one-to-one function from $u_1$ onto $u_2$, then $\pi_{\mathbf{x}}$ is the natural
isomorphism between $\mathbb{Q}_{\mu, u_1}$ and $\mathbb{Q}_{\mu, u_2}$ induced by $h$, i.e. if
$p \in \mathbb{Q}_{\mu, u_1}$ then $\pi_{\mathbf{x}} (p)$ is $q$ with ${\rm dom} (q) = h [{\rm dom} (p)]$
and $q (h(\alpha)) = p (\alpha)$ for $\alpha \in {\rm dom} (p)$.
\end{defn}

\begin{cla}\label{ara}
If $u \subseteq \mu$, then $\mathbb{Q}_{\mu, u}$ is a complete subforcing of $\mathbb{Q}_\mu,\,
\mathbb{Q}_{\mu, u}\hspace{-0.4ex} \leqslant\hspace{-0.7em}\raisebox{0.2ex}{$\cdot$}\, \mathbb{Q}_\mu$ for short. More exactly,
if $p \in \mathbb{Q}_\mu$ then
\begin{itemize}
\item[(a)] if $p \upharpoonright u := p \upharpoonright (u \cap {\rm dom} (p))$, then
$p \upharpoonright u \in \mathbb{Q}_{\mu, u}$ and $p \upharpoonright u \leqslant_{\mathbb{Q}_\mu} p$;
\item[(b)] if $q \in \mathbb{Q}_{\mu, u}$ and $p \upharpoonright u \leqslant q$, then $p$ and
$q$ are compatible in $\mathbb{Q}_\mu$.
\end{itemize}
\end{cla}

\textbf{Proof of Claim \ref{ara}}
(a) is clear. For (b), we have $n_p \leqslant n_q$. Define $r \in \mathbb{Q}_\mu$ with $n_r = n_q,
{\rm dom} (r) = {\rm dom} (p) \cup {\rm dom} (q), r \upharpoonright u = q$ so that for every $\alpha \in {\rm dom}
(p) \setminus u, p (\alpha) \trianglelefteq r (\alpha)$ and $r (\alpha) (n) = 0$ for every $n \in [n_p,
n_q)$. Then $p \leqslant r$ and $q \leqslant r$ hold.\hfill
$\Box$

\textbf{Remark}
This implies that for every filter $G$ that is $\mathbb{Q}_\mu$-generic over some model,
$G \upharpoonright \mathbb{Q}_{\mu, u} := \{p \upharpoonright u : p \in G\}$ in
$\mathbb{Q}_{\mu, u}$-generic.

\begin{defn}\label{bison}
For $\mathbb{C}$ a p.m.s. we define $\mathbb{Q} = \mathbb{Q}_{\mathbb{C}}$ as the finite
support product of $\langle \mathbb{Q}_\mu : \mu \in \mathbb{C}\rangle$.
\end{defn}

Forcing $\mathbb{Q}$ has many natural complete subforcings. In order to talk about them
we introduce the following notations:

\begin{defn}\label{emu}
Let $\mathbb{C}$ be a p.m.s.

\begin{itemize}
\item[1)] For $C \subseteq \mathbb{C}$ we let ${\rm par}_C = \{\overline{u} : \overline{u}
= \langle u_\mu : \mu \in C\rangle$ and $\forall \mu \in C\, \, u_\mu \subseteq \mu\}$ and
then ${\rm par}_{\mathbb{C}} = \bigcup \{ {\rm par}_C : C \subseteq \mathbb{C}\}$
\item[2)] For $\overline{u} \in {\rm par}_{\mathbb{C}}$ let $\mathbb{Q}_{\overline{u}} =
\mathbb{Q}_{\mathbb{C}, \overline{u}}$ be $\mathbb{Q}_\mathbb{C}$ restricted to $\{p \in
\mathbb{Q}_{\mathbb{C}} : {\rm dom}\, (p) \subseteq {\rm dom}\, (\overline{u})$ and
$\forall \mu \in {\rm dom}\, (p)\, p (\mu) \in \mathbb{Q}_{\mu, u_\mu}\}$.
\item[3)] For $\overline{u} \in {\rm par}_{\mathbb{C}}$ and $p \in \mathbb{Q}_{\mathbb{C}}$ let
$p \upharpoonright \overline{u}$ be $q \in \mathbb{Q}_{\overline{u}}$ defined by: ${\rm dom}\, (q)
= {\rm dom}\, (p) \cap {\rm dom}\, (\overline{u})$ and $\forall \mu \in {\rm dom}\, (q)\, q (\mu) =
p (\mu) \upharpoonright u_\mu \cap {\rm dom}\, (p (\mu))$.
\item[4)] We consider partial automorphisms of $\mathbb{Q}_\mathbb{C}$, i.e. ones between
subforcings of the form $\mathbb{Q}_{\overline{u}}$ for $\overline{u} \in {\rm par}_\mathbb{C}$.
We let ${\rm paut}_\mathbb{C}$ be the set of all $\mathbf{x}$ of the form $\langle g, \overline{h},
C_1, \overline{u}_1, C_2, \overline{u}_2\rangle = \langle g_\mathbf{x}, \overline{h}_\mathbf{x},
C_{\mathbf{x}, 1}, \overline{u}_{\mathbf{x}, 1}, C_{\mathbf{x}, 2}, \overline{u}_{\mathbf{x}, 2}\rangle$
such that

\begin{itemize}
\item[(a)] $C_1, C_2, \subseteq \mathbb{C}$;
\item[(b)] $g$ is a one-to-one function from $C_1$ onto $C_2$;
\item[(c)] $\overline{u}_l = \langle u_{l, \mu} : \mu \in C_l \rangle \in {\rm par}_{C_l}$
for $l = 1, 2$;
\item[(d)] $\overline{h} = \langle h_\mu : \mu \in C_1 \rangle$;
\item[(e)] if $g (\mu_1) = \mu_2$, then $h_{\mu_1}$ is a one-to-one function from $u_{1, \mu_1}$
onto $u_{2,\mu_2}$.
\end{itemize}

\item[5)] For $\mathbf{x} \in {\rm paut}_\mathbb{C}$ let $\kappa_\mathbf{x}$ be the isomorphism
between $\mathbb{Q}_{\overline{u}_{\mathbf{x}, 1}}$ and $\mathbb{Q}_{\overline{u}_{\mathbf{x}, 2}}$
which is induced by $\mathbf{x}$.
\end{itemize}
\end{defn}

Generalizing claim \ref{ara}, we easily see that $\mathbb{Q}_{\overline{u}}$
is a complete subforcing of $\mathbb{Q}_\mathbb{C}$:

\begin{cla}\label{frida}
If $\overline{u} \in {\rm par}_{\mathbb{C}}$ then $\mathbb{Q}_{\overline{u}} \hspace{-0.4ex} \leqslant\hspace{-0.7em}\raisebox{0.2ex}{$\cdot$}\,
\mathbb{Q}_\mathbb{C}$. More exactly, if $p \in \mathbb{Q}_{\mathbb{C}}, q \in \mathbb{Q}_{\overline{u}}$
and $\mathbb{Q}_{\overline{u}} \models p \upharpoonright \overline{u} \leqslant q$, then $p$ and $q$
are compatible in $\mathbb{Q}_{\mathbb{C}}$.
\end{cla}

\begin{defn}\label{jojo}
1) Let $\utilde{\mathbb{G}}_{\mathbb{Q}_\mu}$ be the canonical name for the $\mathbb{Q}_\mu$-generic
filter, and let $\utilde\eta_{\mu, \alpha}$ be the $\mathbb{Q}_\mu$-name $\bigcup \{p (\alpha) :
p \in \utilde{\mathbb{G}}_{\mathbb{Q}_\mu}\}$.\\
2) For $\alpha < \mu$ let $\utilde{A}_{\mu, \alpha}$ be the $\mathbb{Q}_\mu$-name $\{n : \utilde{\eta}_{\mu, \alpha}
(n) = 1\}$ and $\utilde{A}_\mu = \langle \utilde{A}_{\mu, \alpha} : \alpha < \mu\rangle$.\\
3) We can consider all these names as $\mathbb{Q}_C$-names, or as $\mathbb{Q}_{\overline{u}}$-names,
provided that $\mu \in C$ or $\overline{u} \in {\rm par}_\mathbb{C}, \mu \in {\rm dom}\, (\overline{u})$
and $\alpha \in u_\mu$ respectively.
\end{defn}

\begin{prop}\label{dicht}
\begin{itemize}
\item[(1)] $\mathbb{Q}_\mathbb{C}$ has the Knaster property and is of cardinality $\max (\mathbb{C})$ such that
$\Vdash_{\mathbb{Q}_{\mathbb{C}}} 2^{\aleph_0} = \max (\mathbb{C})$.
\item[(2)] $\Vdash_{\mathbb{Q}_\mu}\, '' \utilde{\eta}_{\mu, \alpha} \in {^\omega} 2$ and $\utilde{A}_\mu$
is a mad family on $\omega''$, for every $\alpha < \mu$.
\item[(2 A)] $\Vdash_{\mathbb{Q}_{\mathbb{C}}}\, '' \utilde{\eta}_{\mu, \alpha} \in {^\omega} 2$ and
$\utilde{A}_\mu$ is a mad family$''$, for every $\mu \in \mathbb{C}$ and $\alpha < \mu$.
\item[(3)] If $\Vdash_{\mathbb{Q}_\mu} \utilde{\nu} \in {^\omega} 2$, then there are $\alpha_n < \mu
(n < \omega)$ and a Borel function $\mathbb{B} : {^\omega} ({^\omega} 2) \rightarrow {^\omega} 2$,
such that $\Vdash_{\mathbb{Q}_\mu} \utilde{\nu} = \mathbb{B} (\utilde{\eta}_{\alpha_0}, \dots,
\utilde{\eta}_{\alpha_n}, \dots)$.
\item[(3 A)] If $\utilde{\chi}$ is a $\mathbb{Q}_\mu$-name for a subset of an ordinal $\gamma$, then
there is $u \subseteq \mu$ such that $\utilde{\chi}$ is a $\mathbb{Q}_{\mu, u}$-name and $|u| \leqslant
| \gamma | + \aleph_0$.
\item[(4)] If $\Vdash_{\mathbb{Q}_\mathbb{C}}\, \utilde{\nu} \in {^\omega} 2$, then there are
$\mu_n \in \mathbb{C}\, \, \alpha_n < \mu_n (n < \omega)$ and a Borel function $\mathbb{B} : {^\omega}
({^\omega} 2) \rightarrow {^\omega} 2$ such that $\Vdash_{\mathbb{Q}_\mathbb{C}} \utilde{\nu} =
\mathbb{B} (\utilde{\eta}_{\mu_0, \alpha_0}, \dots, \utilde{\eta}_{\mu_n, \alpha_n}, \dots)$.
\item[(4 A)] If $\utilde{\chi}$ is a $\mathbb{Q}_{\mathbb{C}}$-name for a subset of some ordinal
$\gamma$, then for some ordinal $\overline{u} \in {\rm par}_{\mathbb{C}}, \utilde{\chi}$ is a
$\mathbb{Q}_{\overline{u}}$-name and $\mathop{\Sigma}\limits_{\mu \in {\rm dom}\, \overline{u}} | u_\mu |
+ 1 \leqslant | \gamma | + \aleph_0$.
\end{itemize}
\end{prop}

\textbf{Proof:} All arguments needed form part of the basic theory of forcing. Therefore we only
give some hints.
\begin{itemize}
\item[(1)] The Knaster property is preserved by finite support products. By $|\mathbb{Q}_{\mathbb{C}}| = \max
(\mathbb{C})^{\aleph_0} = \max (\mathbb{C})$ (see 1.1 (d)) and the c.c.c. we conclude
$\Vdash_{\mathbb{Q}_\mathbb{C}} 2^{\aleph_0} \leqslant \max (\mathbb{C})$. The converse
follows from (2) below.

\item[(2)] The proof that in (2) and (2 A) $\utilde{A}_\mu$ is forced to be an a.d. family is an
easy genericity argument. Let us prove maximality. Suppose that $p \in \mathbb{Q}_\mathbb{C}$
and $\utilde{a}$ is a $\mathbb{Q}_{\mathbb{C}}$-name such that $p \Vdash_{\mathbb{Q}_\mathbb{C}}
{''} \utilde{a} \in [\omega]^\omega$ and $\utilde{a} \notin \utilde{A}_\mu$ and $\utilde{A}_\mu
\cup \{\utilde{a}\}$ is a. d. $''$. By the c.c.c. of $\mathbb{Q}_\mathbb{C}$ we can find
$\overline{u} \in {\rm par}_\mathbb{C}$ such that $\mathop{\Sigma}\limits_{\nu \in {\rm dom}\,
(\overline{u})} | u_\nu | + 1 \leqslant \aleph_0,\, p \in \mathbb{Q}_{\overline{u}}$ and
$\utilde{a}$ is a $\mathbb{Q}_{\overline{u}}$-name. Fix $\alpha \in \mu \setminus u_\mu$ and find
$q \in \mathbb{Q}_{\mathbb{C}}$ and $m < \omega$ such that $q \geqslant p$ and $q
\Vdash_{\mathbb{P}_\mathbb{C}}\, \utilde{a} \cap \utilde{a}_{\mu, \alpha} \subseteq m$.
By our assumptions we can choose $k \geqslant m$ and $p_1 \in \mathbb{Q}_{\overline{u}}$ such
that $k \geqslant n_{q (\mu)},\, p_1 \geqslant q \upharpoonright \overline{u}$ and $p_1
\Vdash_{\mathbb{Q}_{\overline{u}}} {''} k \in \utilde{a} \setminus \utilde{a}_{\mu, \beta}{''}$
for all $\beta \in u_\mu \cap {\rm dom}\, (q (\mu))$. Note that then $n_{p_1 (\mu)} > k$. We
define $q_1 \in \mathbb{Q}_\mathbb{C}$ as follows: $q_1 \upharpoonright \mathbb{Q}_{\overline{u}}
= p_1,\, q_1 (\nu) = q (\nu)$ for all $\nu \in {\rm dom}\, (q) \setminus {\rm dom}\, (p_1)$
($\mu \in {\rm dom}\, (p_1)$ clearly), $q_1 (\mu) (\beta) \upharpoonright n_{q (\mu)} = q (\mu)
(\beta)$ for all $\beta \in {\rm dom}\, (q (\mu)), q_1 (\mu) (\beta) (n) = 0$ for all $\beta \in {\rm dom}\,
(q (\mu)) \setminus ({\rm dom}\, (p_1 (\mu)) \cup \{\alpha\})$ and $n \in [n_{q (\mu)}, n_{p_1 (\mu)})
$ and finally (the crucial point) $q_1 (\mu) (\alpha) (k) = 1$ and $q_1 (\mu) (\alpha) (n) = 0$ for
all $n \in [n_{q (\mu)}, n_{p_1 (\mu)}) \setminus \{k\}$. Note that $q_1 \geqslant q,\, q_1\geqslant
p_1$ and $q_1 \Vdash_{\mathbb{Q}_\mathbb{C}}\, k \in \utilde{a} \cap \utilde{a}_{\mu, \alpha}$.
This contradicts our choice of $q, m$ and $k$. This finishes the proof of (2) and (2 A).

\item[(3)] We can choose maximal antichains $A_n \subseteq \mathbb{Q}_\mu$ and functions
$f_n : A_n \rightarrow {^{n+1}} 2\, (n < \omega)$ such that $A_{n+1}$ refines $A_n$ and $
\forall_n \forall\, p \in A_n\, p \Vdash_{\mathbb{Q}_\mu} \utilde{\nu} \upharpoonright
n+1 = f_n (p)$. Let $u = \bigcup \{{\rm dom}\, (p) : p \in A_n, n < \omega\}$. By the
c.c.c. we have $|u| \leqslant \aleph_0$. We can consider $\utilde{\nu}$ as a
$\mathbb{Q}_{\mu, u}$-name, and for every $\mathbb{Q}_\mu$-generic filter $G$ we have
$\utilde{\nu} [G] = \utilde{\nu} [G\cap \mathbb{Q}_{\mu, u}]$. Each $p \in \mathbb{Q}_{\mu, u}$
obviously determines a basic open set $U_p$ in the product topology on ${^u} ({^\omega} 2)$.

By the remark after Definition \ref{hei} we need not have $U_p \cap U_q = \emptyset$ for
$p, q \in A_n, p \neq q$. That is why for $n < \omega$ and $p \in A_n$ we let $V_p = U_p
\setminus \bigcup \{U_q : q \in A_n\, \, {\&}\, \, p \neq q\}$. Clearly, $V_p$ is $G_\delta$,
and $V_p \cap V_q = \emptyset$ for any distinct $p, q \in A_n$.

If we let $W_n = \bigcup \{V_p : p \in A_n\}$, we have $W_{n+1} \subseteq W_n$ for all
$n, \bigcap\limits_{n < \omega} W_n$ is $G_{\delta \sigma \delta}$, and for every
$\overline{x} \in \bigcap\limits_{n < \omega} W_n$ and $n < \omega$ there exists a unique
$p \in A_n$ with $\overline{x} \in V_p$. Therefore the functions $f_n$ induce a function
$B' : \bigcap\limits_{n < \omega} W_n \rightarrow {^\omega}2$. Note that its preimage
of any basic open set in ${^\omega} 2$ is $G_{\delta \sigma}$. Hence, if we define
$\mathbb{B}'$ to be constantly zero on ${^\omega} 2  \setminus \bigcap\limits_{n < \omega}
W_n$, then $B' : {^u}({^\omega} 2) \rightarrow {^\omega} 2$ is Borel.

If $\langle \alpha_n : n < \omega\rangle$ is an
enumeration of $u$ and $g : {^\omega} ({^\omega} 2) \rightarrow {^u} ({^\omega} 2), (x_n)
\mapsto (y_\alpha)$ where $y_{\alpha_n} = x_n$, then $\mathbb{B} : = \mathbb{B}' \circ g$
is the desired Borel function. The remaining clauses can be proved by arguments similar
to the ones we used so far.\hfill
$\Box$
\end{itemize}

\newpage

{\sect {\bf Eliminating  $\bigcup \mathbb{C} \cap \mathbb{C}{\rm ard} \subseteq \mathbb{C}$}}

\begin{theo}\label{golo}
Suppose that $\mathbb{C}$ is a p.m.s. such that

\begin{itemize}
\item[(a)] $\min (\mathbb{C}) = \aleph_1$ and $2^{\aleph_0} = \aleph_1$, and
\item[(b)] $\max (\mathbb{C})^{\aleph_0} = \max (\mathbb{C})$.
\end{itemize}

There exists a forcing $\mathbb{Q}_\mathbb{C}$ with the Knaster condition such that,
letting $\utilde{\mathcal{A}}$ a $\mathbb{Q}_\mathbb{C}$-name for the mad spectrum
in $\mathbf{V}^{\mathbb{Q}_\mathbb{C}}$, we have $\Vdash_{\mathbb{Q}_\mathbb{C}}
\utilde{\mathcal{A}} = \mathbb{C}$.
\end{theo}

\textbf{Proof:} Let $\mathbb{Q} = \mathbb{Q}_\mathbb{C}$ as in Definition \ref{bison}.
By Proposition \ref{dicht} (1) we have $\Vdash_\mathbb{Q} 2^{\aleph_0} = \max (\mathbb{C})$.
Let $\lambda \notin \mathbb{C}, \lambda < \max (\mathbb{C})$ be an uncountable cardinal
in $\mathbf{V}^\mathbb{Q}$, hence by the c.c.c. of $\mathbb{Q}$ also in $\mathbf{V}$.
By property (e) of a p.m.s. there exists a minimal regular uncountable cardinal
$\sigma \leqslant \lambda$ such that $[\sigma, \lambda] \cap \mathbb{C} = \emptyset$.
Letting $\chi = \min \{\mu : \mu^{\aleph_0} \geqslant \sigma\}$, we have either

\textbf{Case A:} $\chi = \sigma$, or\vspace{1ex}\\
\textbf{Case B:} $\chi > \aleph_1$ and $cf (\chi) = \aleph_0$.

Indeed, if $\chi < \sigma$, then certainly $\chi > \aleph_1$, as $\aleph_1^{\aleph_0}
= \aleph_1 \in \mathbb{C}$ by assumption. If we had $cf (\chi) \geqslant \aleph_1$,
then $\chi^{\aleph_0} = \mathop{\Sigma}\limits_{\alpha < \chi} |\alpha|^{\aleph_0} < \sigma$, a
contradiction to the definition of $\chi$.

We shall prove Case B and then indicate how the proof can be simplified to
treat Case A.

Assume $p \Vdash_\mathbb{Q} {''} \langle \utilde{B}_\alpha : \alpha < \lambda\rangle$
is an a.d. family${''}$. We have to define a $\mathbb{Q}$-name $\utilde{B}_\lambda$ so that for
every $\alpha < \lambda$

\[p \Vdash_\mathbb{Q}\, \,  ''\utilde{B}_\lambda \in [\omega]^\omega\, \, {\rm and\, \, }
\utilde{B}_\alpha, \utilde{B}_\lambda\, \, {\rm are\, \,  a.d.}''.\]

For this we shall construct $\utilde{B}_\lambda$ with the property that for every
$\alpha < \lambda$ we can find $\beta \in \sigma \setminus \{\alpha\}$ and
$\mathbf{y} = \langle g, \overline{h}, C_1, \overline{u}_1, C_2, \overline{u}_2\rangle
\in {\rm paut}_\mathbb{C}$ (see Definition \ref{bison} (4)) so that

\begin{itemize}
\item[$(\ast)_1$]
\begin{itemize}
\item[(a)] $\utilde{B}_\alpha, \utilde{B}_\beta$ are $\mathbb{Q}_{\overline{u}_1}$-names;
\item[(b)] $\utilde{B}_\alpha, \utilde{B}_\lambda$ are $\mathbb{Q}_{\overline{u}_2}$-names;
\item[(c)] $p \in \mathbb{Q}_{\overline{u}_1} \cap \mathbb{Q}_{\overline{u}_2}$ and
$\kappa_\mathbf{y} (p) = p$;
\item[(d)] $\kappa_\mathbf{y}$ maps $\utilde{B}_\alpha$ to $\utilde{B}_\alpha$ and
$\utilde{B}_\beta$ to $\utilde{B}_\lambda$.
\end{itemize}
\end{itemize}

Since $\kappa_\mathbf{y}$ respects the forcing relation, this will suffice. We have to find
$\utilde{B}_\lambda$ as desired.

By applying Proposition \ref{dicht} (4), for every $\alpha < \lambda$
we can find $\mu (\alpha, n) \in \mathbb{C},\, \xi (\alpha, n, m) < \mu (\alpha, n)$
for $n, m < \omega$ and Borel functions $\mathbb{B}_\alpha$ such that $\Vdash_\mathbb{Q}
{''}\utilde{B}_\alpha = \mathbb{B}_\alpha (\dots, \utilde{\eta}_{\mu (\alpha, n), \xi (\alpha, n, m)},
\dots)_{n, m}{''}$.

For notational simplicity we may assume that all families $\langle \mu (\alpha, n) : n < \omega\rangle
(\alpha < \lambda)$ and $\langle \xi (\alpha, n, m) : m < \omega\rangle (\alpha < \lambda, n < \omega)$
are with no repetition. For each $\alpha < \lambda$ we assemble these ordinals into one sequence
$\overline{\zeta}_\alpha = \langle \zeta (\alpha, \nu) : \nu < \omega \cdot \omega \rangle$ by letting
$\zeta (\alpha, n) = \mu (\alpha, n)$ for $n < \omega$ and $\zeta (\alpha, \omega \cdot (n+1) + m) =
\xi (\alpha, n, m)$ for $n, m < \omega$.

We claim that we can find an unbounded set $Y \subseteq \sigma$, a Borel function $\mathbb{B}_*$, a
partition $\omega \cdot \omega = w_0\, \dot\cup\, w_1\,  \dot\cup\, w_2$ and an ordinal function
$\overline{\beta}^* = \langle \beta^* (i) : i \in w_0\, \cup\, w_1\rangle$ such that

\begin{itemize}
\item[$(\ast)_2$]
\begin{itemize}
\item[(a)] $\mathbb{B}_\alpha = \mathbb{B}_*$ for every $\alpha \in Y$;
\item[(b)] $cf (\beta^* (i)) > \aleph_1$ for every $i \in w_1$;
\item[(c)] for every $\alpha \in Y$
\begin{itemize}
\item[$(\alpha)$] $\overline{\zeta}_\alpha \upharpoonright w_0 = \overline{\beta}^* \upharpoonright w_0$,
\item[$(\beta)$] $\zeta (\alpha, i) < \beta^*_i$ for every $i \in w_1$;
\end{itemize}
\item[(d)] if $\overline{\gamma} \in \prod \overline{\beta}^* \upharpoonright w_1$, then for
$\sigma$ many $\alpha < \sigma$ we have $\overline{\gamma} < \overline{\zeta}_\alpha \upharpoonright
w_1$ (i.e. $\gamma_i < \zeta (\alpha, i)$ for every $i \in w_1$);
\item[(e)] for every $\alpha \in Y$ and $i \in w_2$\\
 $\zeta (\alpha, i) \notin \{\zeta (\beta, j) : \beta < \alpha, j < \omega\cdot\omega\}$.
\end{itemize}
\end{itemize}

In Case A we shall have $w_1 = \emptyset$, hence only (a), (c) ($\alpha$) and (e) are relevant.
We prove $(*)_2$: As there are only $2^{\aleph_0}$ Borel functions and we assume $2^{\aleph_0} =
\aleph_1 < \sigma$, without loss of generality we may assume that $\mathbb{B}_\alpha = \mathbb{B}_*$
for some $\mathbb{B}_*$ and every $\alpha < \sigma$. Let $Z_\alpha = \{\zeta (\beta, \nu) : \beta
< \alpha, \nu < \omega\cdot \omega\}$ and define a function $h$ on $\sigma$ by $h (\alpha) =
\min \{\beta \leqslant \alpha : \forall \nu < \omega \cdot \omega\, (\zeta (\alpha, \nu) \in Z_\alpha
\Rightarrow \zeta (\alpha, \nu) \in Z_\beta)\}$ and let $v_\alpha = \{\nu < \omega\cdot\omega : \zeta (\alpha,
\nu) \in Z_\alpha\}$.

Clearly $h (\alpha) < \alpha$ for $\alpha$ of uncountable cofinality. By Fodor's Lemma there
exist a stationary $S_0 \subseteq \sigma$ and $\gamma < \sigma$ such that $h \upharpoonright S_0$ is
constant with value $\gamma$. Since there are only $\aleph_1$ many possibilities for $v_\alpha$ and
$\sigma$ is regular, there exist a stationary $S_1 \subseteq S_0$ and $v_* \subseteq \omega
\cdot \omega$ such that $v_\alpha = v_*$ for every $\alpha \in S_1$. We let $w_2 = \omega \cdot
\omega \setminus v_*$. Then clearly (e) holds with $S_1$ in place of $Y$. As for $\alpha \in S_1$
and $\nu \in v_*$ we have $\zeta (\alpha, \nu) \in Z_\gamma$ and $|Z_\gamma | \leqslant |\gamma|
\cdot \aleph_0 < \sigma$, in Case A we have $|Z_\gamma|^{\aleph_0} < \sigma$ and hence we can
let $w_0 = v_*$ and find $\overline{\beta}^* = \langle\beta^* (i) : i \in w_0\rangle$ and stationary
$Y \subseteq S_1$ such that (c)($\alpha$) holds.

However, in Case B it may be impossible to make $\overline{\zeta}_\alpha \upharpoonright v_*$
constant for $\sigma$ many $\alpha$, as possibly $|Z_\gamma|^{\aleph_0} \geqslant \sigma$. In
this case we can apply [Sh620, 7.1 (0), (1)] in a straightforward manner with $\lambda, \kappa,
\mu, \mathcal{D}, \langle f_\alpha : \alpha < \lambda\rangle$ there standing for our $S_1, v_*,
\aleph_2,\mathcal{D}_\sigma^{cb}, \langle \overline{\zeta} \upharpoonright v_* : \alpha \in
S_1 \rangle$, where $\mathcal{D}_\sigma^{cb}$ is the filter generated by all cobounded subsets
of $\sigma$. This gives us $Y \subseteq S_1, v = w_0\, \dot\cup\, w_1$ and $\overline{\beta}^*$
as desired.

We are now ready to define the $\mathbb{Q}$-name $\utilde{B}_\lambda$ as outlined at the beginning
of this proof, so that $(*)_1$ will hold. We do it in the Case B, which includes Case A by deleting
everything which refers to $i \in w_1$. We shall define

\[\utilde{B}_\lambda = \mathbb{B}_* (\ldots, \utilde{\eta}_{\mu (\lambda, n), \xi (\lambda, n, m)},
\ldots)_{n, m}\]

for certain $\mu (\lambda, n), \xi (\lambda, n, m)$ which are defined as follows:

\begin{itemize}
\item[$(\ast)_3$]
\begin{itemize}
\item[(a)] If $n \in w_0$, then $\mu (\lambda, n) = \beta^*_n (= \mu (\alpha, n)$ for every $\alpha \in Y)$;
\item[(b)] if $n \notin w_0$, then $\mu (\lambda, n) = \Upsilon_{i (\lambda, n)}$,
where $i (\lambda, n)$ is the $n$-th member of $\{i < \max (\mathbb{C}) : \Upsilon_i \notin \{ \mu (\alpha, k)
: \alpha < \lambda, k < \omega\}\}$;
\item[(c)] if $\omega \cdot (n+1) + m \in w_0$, then $\xi (\lambda, n, m) = \beta^*_{\omega\cdot (n+1)+m} (= \xi (\alpha, m, n)$
for every $\alpha \in Y)$;
\item[(d)] if $n \in w_0$ and $\omega\cdot (n+1) + m \notin w_0$, then $\xi (\lambda, n, m)$
is the $m$-th member of $\mu (\lambda, n) \setminus \{\xi (\beta, n_1, m_1) : \beta < \lambda, n_1, m_1 < w\}$
(Note that this choice is possible, as by $(*)_2$ (d), (e) in this case we must have $\mu (\lambda, n) \geqslant \sigma$
and hence $\mu (\lambda, n) > \lambda$.);
\item[(e)] if $n \notin w_0$ and $\omega \cdot (n+1) + m \notin w_0$, then $\xi (\lambda, n, m) = m$.
\end{itemize}
\end{itemize}

Let $\alpha < \lambda$ be arbitrary. By $(*)_2$ we can choose $\beta \in Y \setminus \{\alpha\}$ such that

\begin{itemize}
\item[$(\ast)_4$]
\begin{itemize}
\item[(a)] if $n \notin w_0$, then $\mu (\beta, n) \notin \{\mu (\alpha, k)
: k < \omega\} \cup {\rm dom}\, (p)$;
\item[(b)] if $\omega \cdot (n+1) + m \notin w_0$, then $\xi (\beta, n, m)
\notin \{\xi (\alpha, n_1, m_1) : n_1, m_1 < \omega\} \cup \bigcup \{ {\rm dom}\,
(p (\mu)) : \mu \in {\rm dom}\, (p)\}$.
\end{itemize}
\end{itemize}

Now we are going to define a partial automorphism of $\mathbb{Q}\, \, \mathbf{y} \in
{\rm paut}_\mathbb{C}$ so that $(*)_1$ will hold.

\newpage

Let the function $g$ be defined by

\begin{itemize}
\item[$(\ast)_5$]
\begin{itemize}
\item[(a)] ${\rm dom}\, (g) = \{\mu (\alpha, n) : n < \omega\} \cup \{\mu
(\beta, n) : n < \omega\} \cup {\rm dom}\, (p)$;
\item[(b)] $g (\mu (\alpha, n)) = \mu (\alpha, n)$;
\item[(c)] $g (\mu (\beta, n)) = \mu (\lambda, n)$, thus\vspace{1ex}
\begin{itemize}
\item[$(\alpha)$] $g (\mu (\beta, n)) = \mu (\beta, n)$ if $n \in w_0$,\vspace{0.5ex}
\item[$(\beta)$] $g (\mu (\beta, n)) = \Upsilon_{i (\lambda, n)}$ if $n \notin w_0$;\vspace{1ex}
\end{itemize}
\item[(d)] $g (\mu) = \mu$ for $\mu \in {\rm dom}\, (p)$.
\end{itemize}
\end{itemize}

Note that by the choice of $\beta$ in $(\ast)_4$, $g$ is well-defined, i.e. if $\mu (\alpha,
n_1) = \mu (\beta, n_2)$ then the demands in (b) and (c) agree, similarly for $\mu (\beta, n_2)
= \mu \in {\rm dom}\, (p)$. Indeed, in this case we must have $n_2 \in w_0$ and $g (\mu (\beta, n_2))
= \mu (\beta, n_2)$.

Let $C_1 := {\rm dom}\, (g)$ and $C_2 := \, {\rm ran}\, (g)$. By definition and by $(*)_3$ (b), $g$
is one-to-one.

For each $\mu \in {\rm dom}\, (g)$ we define a function $h_\mu$ as follows:

\begin{itemize}
\item[$(\ast)_6$]
\begin{itemize}
\item[(a)]
\begin{itemize}
\item[$(\alpha)$] If $\mu = \mu (\alpha, n) \notin \{\mu (\beta, m)
: m < \omega\}$, then ${\rm dom}\, (h_\mu) = \{\xi (\alpha, n, m) : m < \omega\}
\cup {\rm dom}\, (p (\mu))\, \, ({\rm dom}\, (p (\mu)) = \emptyset$ if $\mu \notin {\rm dom} (p))$;\vspace{1ex}
\item[$(\beta)$] if $\mu = \mu (\beta, n) \notin \{\mu (\alpha, m) : m < \omega\}$,
then ${\rm dom}\, (h_\mu) = \{\xi (\beta, n, m) : m < \omega\} \cup {\rm dom}\,
(p (\mu))$;\vspace{1ex}
\item[$(\gamma)$] if $\mu = \mu (\alpha, n_1) = \mu (\beta, n_2)$, then ${\rm dom}\,
(h_\mu) = \{\xi (\alpha, n_1, m) : m < \omega\} \cup \{\xi (\beta, n_2, m) : m <
\omega\} \cup {\rm dom}\, (p (\mu))$;\vspace{1ex}
\item[$(\delta)$] if $\mu \in {\rm dom}\, (p) \setminus \{\mu (\nu, n) : \nu \in
\{\alpha, \beta\}, n < \omega\}$, then ${\rm dom}\, (h_\mu) = {\rm dom}\, (p (\mu))$;
\end{itemize}
\item[(b)] if $\mu = \mu (\alpha, n)$, then $h_\mu (\xi (\alpha, n, m)) = \xi
(\alpha, n, m)$;\vspace{1ex}
\item[(c)] if $\mu = \mu (\beta, n)$, then $h_\mu (\xi (\beta, n, m)) = \xi (\lambda,
m, n)$;\vspace{1ex}
\item[(d)] if $\mu \in {\rm dom}\, (p)$ and $\xi\in {\rm dom}\, (p (\mu))$, then $h_\mu
(\xi) = \xi$.
\end{itemize}
\end{itemize}

Again $h_\mu$ is well defined: E.g. if $\mu = \mu (\alpha, n_1) = \mu (\beta, n_2)$
and $\xi (\alpha, n_1, m_1) = \xi (\beta, n_2, m_2)$, as before we must have $n_2 \in w_0$,
but also $\omega \cdot (n_2 + 1) + m_2 \in w_0$ by $(*)_4$ (b), and hence $h_\mu (\xi (\beta,
n_2, m_2)) = \xi (\lambda, n_2, m_2) = \xi (\beta, n_2, m_2) = \xi (\alpha, n_1, m_1)
= h_\mu (\xi (\alpha, n_1, m_1))$. The other cases are similar. Moreover, $h_\mu$ is one-to-one
by $(*)_3$ (d), (e).

Let $u_{1, \mu} := {\rm dom}\, (h_\mu), u_{2, g (\mu)} := {\rm ran}\, (h_\mu),
\overline{h}:= \langle h_\mu : \mu \in C_1\rangle, \overline{u}_1 := \langle u_{1, \mu}
: \mu \in C_1\rangle, \overline{u}_2 := \langle u_{2, \mu} : \mu \in C_2\rangle, \mathbf{y} =
\langle g, \overline{h}, C_1, \overline{u}_1, C_2, \overline{u}_2\rangle$. Then
$\mathbf{y} \in {\rm paut}_\mathbb{C}$.

We conclude that $\kappa_{\mathbf{y}}$ is an isomorphism between $\mathbb{Q}_{\overline{u}_1}$
and $\mathbb{Q}_{\overline{u}_2}$. By construction we can now easily verify $(\ast)_1$.
Hence Theorem \ref{golo} is proved.\hfill
 $\Box$ \vspace{3ex}

{\sect {\bf Eliminating $\mathbf{\aleph_1 \in \mathbb{C}}$}}

In this section we extend the method of \S 2 to construct a forcing $\mathbb{P}_\mathbb{C}$, for
a given p.m.s. $\mathbb{C}$, so that in addition to Theorem \ref{golo} we can show that in
a $\mathbb{P}_\mathbb{C}$-extension the minimum of the madness spectrum equals $\min (\mathbb{C})$.
For this we first force with $\mathbb{Q}_\mathbb{C}$ from \S 2 and then force a weak from of
$MA_{< \vartheta}$, where $\vartheta = \min (\mathbb{C})$, that rules out mad families of size $< \vartheta$.
More precisely, we shall force $MA_{< \vartheta}$ for all forcings with the Knaster condition. Hence
$\mathbb{P}_\mathbb{C}$ will be the limit of a finite support iteration $\langle \mathbb{P}_\alpha,
\utilde{\mathbb{Q}}_\beta : \alpha \leqslant \vartheta, \beta < \vartheta\rangle$, so $\mathbb{P}_\mathbb{C}
= \mathbb{P}_\vartheta$, where $\mathbb{Q}_0 = \mathbb{Q}_\mathbb{C}$ and $\mathbb{P}_\vartheta / \mathbb{Q}_0$
forces $MA_{< \vartheta}$ (Knaster).

Let us recall that by a well-known reflection argument, for forcing $MA_{< \vartheta}$ (Knaster)
it suffices to take care of posets of size $<\vartheta$ only. However, after forcing with
$\mathbb{Q}_0  = \mathbb{Q}_\mathbb{C}$ we have $2^{\aleph_0} = \max (\mathbb{C})$ (for this we
need that $\max (\mathbb{C})$ exists and $\max (\mathbb{C})^{\aleph_0} = \max (\mathbb{C})$)
(see Definition \ref{coco}). Therefore, forcings $\utilde{\mathbb{Q}}_\beta, 1 \leqslant \beta < \vartheta$,
will be finite support products of length $\max (\mathbb{C})$ of forcings of size $< \vartheta$ with
the Knaster property. Note that the forcing that kills a mad family of size $< \vartheta$ has this
property. For all this to work we have to assume $\vartheta^{<\vartheta} = \vartheta$ and ${\rm max}\,
(\mathbb{C})^{<\vartheta} = {\rm max}\, (\mathbb{C})$.

Moreover we want to preserve what we have obtained in Theorem \ref{golo}. In fact we want to be able
to repeat essentially the same arguments using partial automorphisms as in \S 2. For this goal,
simultaneously to defining the iteration we define (many) names for parameters of partial
automorphisms and for complete subforcings of $\mathbb{P}_\alpha, \utilde{\mathbb{Q}}_\beta
(\alpha \leqslant \vartheta, \beta < \vartheta)$. In this way we shall prove the following Theorem:

\begin{theo}\label{isa}
Assume that $\mathbb{C}$ is a p.m.s., $\vartheta := \min (\mathbb{C})$ satisfies $\vartheta = \vartheta^{<\vartheta}$
and $\max (\mathbb{C})^{<\vartheta} = \max (\mathbb{C})$. Then we can find $\mathbb{P}_\mathbb{C}$ such
that, letting $\utilde{\mathcal{A}}$ be a $\mathbb{P}_\mathbb{C}$-name for the mad spectrum
in $\mathbb{V}^{\mathbb{P}_\mathbb{C}}$, we have:

\begin{itemize}

 \item[(a)] $\mathbb{P}_\mathbb{C}$ is a c.c.c. forcing notion of cardinality $\max (\mathbb{C})$;
\item[(b)] $\Vdash_{\mathbb{P}_\mathbb{C}}\, \mathbb{C} = \utilde{\mathcal{A}}$.
\end{itemize}
\end{theo}

\textbf{Proof:} Recursively we construct the following objects:

\begin{itemize}
\item[$\bullet$] Partial orders $\mathbb{P}_\alpha, I_\alpha, J'_\alpha, J_\alpha,
\mathbb{P}_{\alpha, f}$ for $\alpha \leqslant \vartheta, f \in J_\alpha$;
\item[$\bullet$] $\mathbb{P}_\beta$-names for partial orders $\utilde{\mathbb{Q}}_\beta,
\utilde{\mathbb{Q}}_{\beta, \varepsilon}, \utilde{\mathbb{Q}}_{\beta, t}$ for
$\beta < \vartheta, \varepsilon < \max (\mathbb{C}), t \in I_\beta$;
\item[$\bullet$] elements of $J_\alpha\, \, g_{\alpha, \varepsilon}$ for
$\varepsilon < \max (\mathbb{C})$;
\item[$\bullet$] ordinals $\gamma_{\beta, \varepsilon}$ for $\beta < \vartheta,
\varepsilon < \max (\mathbb{C})$;
\item[$\bullet$] names $\utilde{\nu}_{\beta, \varepsilon}$ for subsets of
$\gamma_{\beta, \varepsilon}$, for $\beta < \vartheta, \varepsilon < \max (\mathbb{C})$.
\end{itemize}

The following properties shall be satisfied:

\begin{itemize}
\item[(A)] $\langle \mathbb{P}_\alpha, \utilde{\mathbb{Q}}_\beta : \alpha \leqslant
\vartheta, \beta < \vartheta\rangle$ is a finite support iteration of forcing notions with
the Knaster condition such that $\Vdash_{\mathbb{P}_\alpha} | \utilde{\mathbb{Q}}_\alpha
| = \max (\mathbb{C})$ for every $\alpha < \vartheta$;

\item[(B)] $\mathbb{Q}_0 = \utilde{\mathbb{Q}}_0 = \mathbb{Q}_\mathbb{C}$ (see defintion
\ref{bison});

\item[(C)]
\begin{itemize}
\item[($\alpha$)] $I_0 = \{\overline{u} \in {\rm par}_\mathbb{C} :
\mathop{\Sigma}\limits_{\mu \in {\rm dom}\, (\overline{u})} | u_\mu | +1 < \vartheta\}$ partially
ordered by $\overline{u}_1 \leqslant_{I_0} \overline{u}_2$ iff ${\rm dom}\, (\overline{u}_1)
\subseteq {\rm dom}\, (\overline{u}_2)$ and $u_{1, \mu} \subseteq u_{2, \mu}$ for all
$\mu \in {\rm dom}\, (\overline{u}_1)$ (see Definition \ref{emu}.1) $I_{1+\alpha} = [\max
(\mathbb{C})]^{< \vartheta}$ from $\mathbf{V}$ partially ordered by $\subseteq$, for every
$\alpha < \vartheta$;
\item[($\beta$)] $I_\alpha$ is $\vartheta$-directed for every $\alpha$;
\end{itemize}

\item[(D)]
\begin{itemize}
\item[($\alpha$)] $J'_\alpha = \{f : {\rm dom}\, (f) \in ([\alpha]^{<\vartheta})^\mathbf{V}$
and if $\beta \in {\rm dom}\, (f)$, then $f (\beta) \in I_\beta\}$;
\item[($\beta$)] $f \leqslant_{J'_\alpha} g$ iff ${\rm dom}\, (f) \subseteq {\rm dom}\, (g)$
and if $\beta \in {\rm dom}\, (f)$, then $f (\alpha) \leqslant_{I_\alpha} g (\alpha)$;
\item[($\gamma$)] $J_\alpha =\{f\in J'_\alpha :$ if $\beta \in {\rm dom}\, (f), \beta
\neq 0$ and $\varepsilon \in f (\beta)$, then $g_{\beta, \varepsilon} \leqslant_{J'_\alpha}
f \upharpoonright \beta\}$ equipped with the induced partial order $\leqslant_{J'_\alpha}
\upharpoonright J_\alpha$;
\item[($\delta$)] $J_\alpha$ is a cofinal subset of $J'_\alpha$;
\item[($\varepsilon$)] if $\beta < \alpha$ then $J_\beta = \{f \in J_\alpha : {\rm dom}\, (f)
\subseteq \beta\}$ and $\leqslant_{J_\beta} = \leqslant_{J_\alpha} \upharpoonright J_\beta$;
\item[($\zeta$)] $J_\alpha$ and $J'_\alpha$ are $\vartheta$-directed partial orders of cardinality
$\max (\mathbb{C})$;
\end{itemize}

\item[(E)]
\begin{itemize}
\item[($\alpha$)] $\utilde{\mathbb{Q}}_{\alpha, t} \hspace{-0.4ex} \leqslant\hspace{-0.7em}\raisebox{0.2ex}{$\cdot$}\, \utilde{\mathbb{Q}}_\alpha$
increase with $t \in I_\alpha$ (i.e. $s <_{I_\alpha} t \Rightarrow\, \Vdash_{\mathbb{P}_\alpha}
\utilde{\mathbb{Q}}_{\alpha, s} \hspace{-0.4ex} \leqslant\hspace{-0.7em}\raisebox{0.2ex}{$\cdot$}\, \utilde{\mathbb{Q}}_{\alpha, t} \linebreak
\hspace{-0.4ex} \leqslant\hspace{-0.7em}\raisebox{0.2ex}{$\cdot$}\, \utilde{\mathbb{Q}}_\alpha)$;
\item[($\beta$)] $\Vdash_{\mathbb{P}_\alpha} \utilde{\mathbb{Q}}_\alpha = \bigcup \{
\utilde{\mathbb{Q}}_{\alpha, t} : t \in I_\alpha\}$,
\end{itemize}

\item[(F)]
\begin{itemize}
\item[($\alpha$)] $\mathbb{P}_{\alpha, f} \hspace{-0.4ex} \leqslant\hspace{-0.7em}\raisebox{0.2ex}{$\cdot$}\, \mathbb{P}_\alpha$ increase with
$f \in J_\alpha$ (i.e. $f <_{J_\alpha} g \Rightarrow \mathbb{P}_{\alpha, f} \hspace{-0.4ex} \leqslant\hspace{-0.7em}\raisebox{0.2ex}{$\cdot$}\,
\mathbb{P}_{\alpha, g} \hspace{-0.4ex} \leqslant\hspace{-0.7em}\raisebox{0.2ex}{$\cdot$}\, \mathbb{P}_\alpha$);
\item[($\beta$)] $\mathbb{P}_\alpha = \bigcup \{\mathbb{P}_{\alpha, f} : f \in J_\alpha\}$;
\item[($\gamma$)] $\langle\mathbb{P}_{\beta, f \upharpoonright \beta},
\utilde{\mathbb{Q}}_{\gamma, f (\gamma)} : \beta \in {\rm dom}\, (f) \cup \{\alpha\}, \gamma
\in {\rm dom}\, (f)\rangle$ is a finite support iteration, for every $f \in J_\alpha$;
\item[($\delta$)] $\mathbb{P}_{\alpha, f}$ has density $< \vartheta$;
\end{itemize}

\item[(G)] $\mathbb{Q}_{0, \overline{u}} = \mathbb{Q}_{\mathbb{C}, \overline{u}}$ for
$\overline{u} \in I_0$ (see Definition \ref{emu}.2);
\item[(H)]
\begin{itemize}
\item[($\alpha$)] for $\alpha > 0, \langle (\utilde{\mathbb{Q}}_{\alpha, \varepsilon},
g_{\alpha, \varepsilon}) : \varepsilon < \max (\mathbb{C})\rangle$ is a sequence of pairs
from $X_\alpha := \{(\utilde{\mathbb{Q}}, f) : f \in J_\alpha, \utilde{\mathbb{Q}}$ is
a $P_{\alpha, f}$-name of a forcing notion satisfying the Knaster condition with set of
elements an ordinal $< \vartheta$ (not only a $\mathbb{P}_\alpha$-name!)$\}$;
\item[($\beta$)] let $\gamma_{\alpha, \varepsilon}$ be the set of elements of
$\utilde{\mathbb{Q}}_{\alpha, \varepsilon}$;
\item[($\gamma$)] each pair from $X_\alpha$ appears $\max (\mathbb{C})$ times in the
sequence from $(\alpha)$;
\item[($\delta$)] $\Vdash_{\mathbb{P}_{\alpha, g_{\alpha, \varepsilon}}}
\utilde{\nu}_{\alpha, \varepsilon} \in {^{\gamma_{\alpha, \varepsilon}}2}$ is
$\utilde{\mathbb{Q}}_{\alpha, \varepsilon}$-generic;
\end{itemize}

\item[(I)]
\begin{itemize}
\item[($\alpha$)] $\Vdash_{\mathbb{P}_\alpha} \utilde{\mathbb{Q}}_\alpha$ is the finite
support product of $\langle\utilde{\mathbb{Q}}_{\alpha, \varepsilon} : \varepsilon
< \max (\mathbb{C})\rangle$;
\item[($\beta$)] for $t \in I_\alpha, \Vdash_{\mathbb{P}_\alpha} \utilde{\mathbb{Q}}_{\alpha, t}$
is the finite support product of $\langle\utilde{\mathbb{Q}}_{\alpha, \varepsilon} : \varepsilon
\in t\rangle$;
\end{itemize}

\item[(J)] letting $\mathbb{P}_\alpha'$ be the set of all $p \in \mathbb{P}_\alpha$ such
that for every $\beta \in {\rm dom}\, (p), p (\beta)$ is an object and not only a $\mathbb{P}_\beta$-name,
and hence $p (\beta)$ is a finite partial function from $\max (\mathbb{C})$ to $\vartheta$ and
$p (\beta) (\varepsilon) < \gamma_{\beta, \varepsilon}$ for $\varepsilon \in {\rm dom}\, (p (\beta))$,
then $\mathbb{P}'_\alpha$ is a dense subset of $\mathbb{P}_\alpha$; similarly, letting
$\mathbb{P}'_{\alpha, f} = \mathbb{P}'_\alpha \cap \mathbb{P}_{\alpha, f}$ for
$f \in J_\alpha, \mathbb{P}'_{\alpha, f}$ is a dense subset of $\mathbb{P}_{\alpha, f}$ of size $< \vartheta$;

\item[(K)] for $\alpha \leqslant \vartheta, f \in J_\alpha$ and $p \in \mathbb{P}'_\alpha$
\begin{itemize}
\item[($\alpha$)] let $p \upharpoonright f$ be defined as the function $q$ such that
${\rm dom}\, (q) = {\rm dom}\, (p) \cap {\rm dom}\, (f)$, and if $\beta \in {\rm dom}\,
(q)$, then $q (\beta) = p (\beta) \upharpoonright f (\beta)$;
\item[($\beta$)] then $p \upharpoonright f \in \mathbb{P}'_{\alpha, f}$ and
$p \upharpoonright f \leqslant_{\mathbb{P}_\alpha} p$;
\item[($\gamma$)] moreover, if $r \in \mathbb{P}'_{\alpha, f}$ and $p \upharpoonright f
\leqslant_{\mathbb{P}_{\alpha, f}} r$, then $p$ and $r$ are compatible in $\mathbb{P}_\alpha$.
\end{itemize}
\end{itemize}

Verifying inductively that this recursion is well-defined and all relevant claims hold
is essentially the same thing as reading and understanding it carefully, thereby using our
assumptions and well-known facts about finite-support iterations. Therefore we shortly sketch
the order of this recursive construction.

The partial orders $I_\alpha$ and $I'_\alpha$ for $\alpha < \vartheta$ are defined directly in (C), (D) ($\alpha$),
($\beta$). Clearly they are all $\vartheta$-directed.

\textbf{Case 1:} $\alpha = 0$.
We just have to define $\mathbb{P}_0$ as the empty forcing notion, $J_0 = \{\emptyset\}\, \utilde{\mathbb{Q}}_{0,
\overline{u}} = \mathbb{Q}_{0, \overline{u}}$ is defined in (G) and $\mathbb{P}_{0, \emptyset} = \mathbb{P}_0$.
All relevant claims can be checked.

\textbf{Case 2:} $\alpha = 1$.
$\mathbb{P}_1$ is defined by (A) and (B), so $\mathbb{P}_1 \cong \mathbb{Q}_0$; we have $J_1 = J'_1$
which is essentially $I_0$. $\mathbb{P}_{1, f}$ for $f \in J_1$ is defined in (F) ($\gamma$), hence
$\mathbb{P}_{1,f} \cong \mathbb{Q}_{0, f(0)} = \mathbb{Q}_{\overline{u}}$ where $\overline{u} = f (0)$.
$\utilde{\mathbb{Q}}_1, \utilde{\mathbb{Q}}_{1,\varepsilon}, \utilde{\mathbb{Q}}_{1,t}$ for $\varepsilon < {\rm max}\,
(\mathbb{C}), t \in I_1$ are defined in (H), (I). Finally $g_{1, \varepsilon}, \gamma_{1, \varepsilon},
\utilde{\nu}_{1, \varepsilon}$ for $\varepsilon < {\rm max}\, (\mathbb{C})$ are defined in (H). All relevant
claims can be checked. Note that $\utilde{\mathbb{Q}}_1$ and $\utilde{\mathbb{Q}}_{1,t}$ are forced to satisfy
the Knaster condition, as the Knaster property is preserved by finite support products.

\textbf{Case 3:} $\mathbf{\alpha}$ \textbf{is a limit ordinal}.
$\mathbb{P}_\alpha$ is defined by (A), $J_\alpha$ is defined by (D) ($\gamma$), it is $\vartheta$-directed
by the induction hypothesis. $\mathbb{P}_{\alpha, f}$ is defined in (F) ($\gamma$) as the limit of a
finite support iteration of forcings with the Knaster condition. $\utilde{\mathbb{Q}}_\alpha, \utilde{\mathbb{Q}}_{\alpha,
\varepsilon}, \mathbb{Q}_{\alpha, t}$ for $\varepsilon < {\rm max}\, (\mathbb{C}), t \in I_\alpha$ are
defined in (H), (I). Finally $g_{\alpha, \varepsilon}, \gamma_{\alpha, \varepsilon}, \utilde{\nu}_{\alpha, \varepsilon}$
are defined in (H). All relevant claims can be checked.

\textbf{Case 4:} $\alpha = \beta + 1.$
All relevant objects are defined in the same order and by the same clauses as in the limit case.

In order to prove Theorem \ref{isa} we shall essentially repeat the arguments from
${\S}$ 2. For this we need a notation for partial isomorphisms of $\mathbb{P}_\vartheta$.
This will extend Definition \ref{emu}.

\begin{defn}
1) For $0 < \alpha \leqslant \delta$ we define ${\rm p paut}_\alpha$ (for preliminary
partial automorphism) to be the set of all $\mathbf{s} = (f_1, f_2, \mathbf{x}, \overline{k})$
such that
\begin{itemize}
\item[(a)] $f_1, f_2 \in J_\alpha$ satisfy $0 \in {\rm dom}\, (f_1) = {\rm dom}\,
(f_2)$;
\item[(b)] $\mathbf{x} \in {\rm paut}_\mathbb{C}, \overline{u}_{\mathbf{x}, 1} = f_1 (0),
\overline{u}_{\mathbf{x}, 2} = f_2 (0)$ (see Definition \ref{emu}.4);
\item[(c)] $\overline{k} = \langle k_\beta : \beta \in {\rm dom}\, (f_1) \setminus
\{0\}\rangle$ and $k_\beta$ is a bijection from $f_1 (\beta)$ onto $f_2 (\beta)$;
\item[(d)] if $\beta \in {\rm dom}\, (f_1), \varepsilon_1 \in f_1 (\beta)$
and $\varepsilon_2 = k_\beta (\varepsilon_1)$ (hence $g_{\beta, \varepsilon_1}
\mathop{\leqslant}_{\mathop{J}\limits_\alpha} f_1$), then $g_{\beta, \varepsilon_1}$ is
mapped to $g_{\beta, \varepsilon_2}$ by $\mathbf{s}$, which means
\begin{itemize}
\item[($\alpha$)] ${\rm dom}\, (g_{\beta, \varepsilon_1}) = {\rm dom}\, (g_{\beta, \varepsilon_2})$,
\item[($\beta$)] if $\gamma \in {\rm dom}\, (g_{\beta, \varepsilon_1})$, then $k_\gamma [g_{\beta, \varepsilon_1}
(\gamma)] = g_{\beta, \varepsilon_2}(\gamma)$.
\end{itemize}
\end{itemize}

Then we write $\bm{s} = (f_{\bm{s}, 1}, f_{\bm{s}, 2}, \mathbf{x}_{\bm{s}},
\overline{k}_{\bm{s}}), k_\beta = k_{{\bm{s}}, \beta}, \overline{u}_{\mathbf{x}_{\bm{s}}, 1}
= \overline{u}_{\bm{s}, 1}, \overline{u}_{\mathbf{x}_{\bm{s}}, 2} = \overline{u}_{\bm{s}, 2}$.

2) For $1 \leqslant \alpha < \beta \leqslant \delta, \bm{s} \in {\rm p paut}_\beta,
\bm{t} \in {\rm p paut}_\alpha$ we define $\bm{t} = \bm{s} \upharpoonright \alpha$ by
$\bm{t} = (f_{\bm{s}, 1} \upharpoonright {\rm dom}\, (f_{\bm{s}, 1}) \cap \alpha,
f_{\bm{s}, 1} \upharpoonright {\rm dom}\, (f_{\bm{s}, 2}) \cap \alpha,
\mathbf{x}_{\bm{s}}, \overline{k}_{\bm{s}} \upharpoonright ({\rm dom}\, (f_{\bm{s}, 1})
\cap \alpha \setminus \{0\})$.

3) For $f \in J_\alpha$ we let $J_{\alpha,f} = \{g \in J_\alpha : g \leqslant_{J_\alpha}
f\}$. For $\bm{s} \in {\rm p paut}_\alpha$ we define an isomorphism $\pi_{\bm{s}}$
from $J_{\alpha, f_{\bm{s}, 1}}$ onto $J_{\alpha, f_{\bm{s}, 2}}$ by letting
$\pi_{\bm{s}} (g_1) = g_2$ iff:
\begin{itemize}
\item[(a)] $g_l \in J_\alpha$ and $g_l \leqslant_{J_\alpha} f_{\bm{s}, l}$ for
$l = 1, 2$;
\item[(b)] ${\rm dom}\, (g_1) = {\rm dom}\, (g_2)$;
\item[(c)] if $0 \in {\rm dom}\, (g_1)$ then $\mathbf{x}_{\bm{s}}$ naturally maps $g_1 (0)$
to $g_2 (0)$, i.e. letting $\overline{h}_\mathbf{x} = \langle h_\mu : \mu \in
{\rm dom}\, (f_{\bm{s}, 1})$ and $g_l (0) = \langle u_{l, \mu} : \mu \in
{\rm dom}\, (f_{\bm{s}, l})\rangle$ for $l = 1, 2$ (see \ref{emu}.4), $h_\mu$ is a
one-to-one map from $u_{1, \mu}$ onto $u_{2, \mu}$;
\item[(d)] if $\beta \in {\rm dom}\, (g_1) \setminus \{0\}$ then $\{k_\beta
(\varepsilon) : \varepsilon \in g_1 (\beta)\} = g_2 (\beta)$.
\end{itemize}

4) For every $\bm{s} \in {\rm p paut}_\alpha$ we can naturally define
$\bm{s}^{-1} \in {\rm p paut}_\alpha$, so that $\pi_{\bm{s}^{-1}} =
(\pi_{\bm{s}})^{-1}$. Note that for $g_1 \leqslant f_{\bm{s}, 1}$,
if $\pi_{\bm{s}} (g_1) = g_2$ then $\bm{s}^{-1} \upharpoonright g_2 =
(\bm{s} \upharpoonright g_1)^{-1}$.
\end{defn}

5) For $\bm{s} \in {\rm paut}_\alpha$ and $g \leqslant f_{\bm{s}, 1}$ we
define $\bm{s} \upharpoonright g \in {\rm p paut}_\alpha$ in the canonical
way.

\begin{defn}\label{logo}
By recursion on $\alpha \in [1, \vartheta]$ we define ${\rm paut}_\alpha \subseteq
{\rm p paut}_\alpha$ such that $\forall \beta \forall \bm{s} (1 \leqslant
\beta < \alpha \wedge \bm{s} \in {\rm paut}_\alpha) \rightarrow \bm{s}
\upharpoonright \beta \in {\rm paut}_\beta$, and for $\bm{s} \in {\rm paut}_\alpha$
we define an isomorphism $\kappa_{\bm{s}}$ from $\mathbb{P}_{\alpha, f_{\bm{s}, 1}}$
onto $\mathbb{P}_{\alpha, f_{\bm{s}, 2}}$ with the property that if $g \in \mathbb{P}_{\alpha,
f_{\bm{s}, 1}}$ then $\kappa_{\bm{s}} \upharpoonright \mathbb{P}_{\alpha, g} = \kappa_{\bm{s}
\upharpoonright g}$. For $\alpha = 1$ we let ${\rm paut}_\alpha
= {\rm p paut}_\alpha$ and $\kappa_{\bm{s}} = \kappa_{\mathbf{x}_{\bm{s}}}$ (see
Definition \ref{emu}.5).

In case $\alpha = \beta + 1$ for some $1 \leqslant \beta
< \vartheta$ we let ${\rm paut}_\alpha$ be the set of all $\bm{s} \in {\rm p paut}_\alpha$
such that $\bm{s} \upharpoonright \beta \in {\rm paut}_\beta$ and if $\mu \in {\rm dom}\,
(f_{\bm{s}, 1}), \varepsilon_1 \in f_{\bm{s}, 1} (\mu)$ and $\varepsilon_2 =
k_\mu (\varepsilon_1)$ (hence $\varepsilon_2 \in f_{\bm{s}, 2} (\mu))$, then $\gamma_{\mu,
\varepsilon_1} = \gamma_{\mu, \varepsilon_2}$ (the domains of
$\utilde{\mathbb{Q}}_{\mu, \varepsilon_1}, \utilde{\mathbb{Q}}_{\mu, \varepsilon_2})$ and the pair
$\kappa_{\bm{s} \upharpoonright g_{\beta, \varepsilon_1}}, id_{\gamma_{\mu, \varepsilon_1}}$
maps $\utilde{\mathbb{Q}}_{\beta, \varepsilon_1}$ onto $\utilde{\mathbb{Q}}_{\beta, \varepsilon_2}$,
i.e. for every $p \in \mathbb{P}_{\beta, g_{\beta, \varepsilon_1}}$ and $\nu_0, \nu_1 < \gamma_{\mu,
\varepsilon_1}$ we have $p\, {\mathop{\Vdash}_{\mathbb{P}_{\beta, g_{\beta, \varepsilon_1}}}}
''\nu_0 \mathop{<}_{\utilde{\mathbb{Q}}_{\beta, \varepsilon_1}} {\nu_1}''$ iff $\kappa_{\bm{s}} (p)\,
{\mathop{\Vdash}_{\mathbb{P}_{\beta, g_\beta, \varepsilon_2}}} ''\nu_0 \mathop{<}_{\utilde{\mathbb{Q}}_{\beta,
\varepsilon_2}} {\nu_1}''$. In case $\alpha \in [1, \vartheta]$ is a limit ordinal let
${\rm paut}_\alpha = \{\bm{s} \in {\rm p paut}_\alpha : \forall \beta \bm{s} \upharpoonright
\beta \in {\rm paut}_\beta$.

\end{defn}

The following claim extends Proposition \ref{dicht}.4A.

\begin{cla}\label{fix}
If $\utilde{B}$ is a $\mathbb{P}_\alpha$-name $(\alpha \leqslant \delta)$
for a bounded subset of $\delta$, then for some $f \in J_\alpha, \utilde{B}$ is
a $\mathbb{P}_{\alpha, f}$-name.
\end{cla}

\textbf{Proof of claim \ref{fix}:} By the ccc of $\mathbb{P}_\alpha$ and our
assumption $\vartheta = \vartheta^{<\vartheta}$ there exists $\gamma < \vartheta$ such that
$\Vdash_{\mathbb{P}_\alpha} \utilde{B} \subseteq \gamma$. Hence $\utilde{B}$ is
determined by a $\gamma$-sequence of maximal antichains of $\mathbb{P}_\alpha$.
By properties $(D \zeta)$ and $(J)$ it suffices to find for given $p \in
\mathbb{P}'_\alpha$ some $f \in J_\alpha$ with $p \in \mathbb{P}_{\alpha, f}$.
This is trivial.\hfill
$\Box_{\ref{logo}}$\vspace{3ex}

\textbf{Remark:}
\emph{Note that even if $\utilde{B}$ is a name for a real it is generally impossible to
obtain a countable $f$ as in \ref{fix}. The reason is our definition of $J_\alpha$
in (D)($\gamma$).}

We are now ready to prove Theorem \ref{isa}.\\
\ref{isa} (a) follows from (A) and our assumptions about $\vartheta (\vartheta = \min (\mathbb{C})$
and $\vartheta = \vartheta^{< \vartheta})$. In order to prove $\mathop{\Vdash}_{\mathbb{P}_\vartheta}
\mathbb{C} \subseteq \utilde{\mathcal{A}}$ we use the notation from Definition
\ref{jojo} to denote the objects added by $\mathbb{Q}_0 = \mathbb{Q}_\mathbb{C}$.
Hence for $\mu \in \mathbb{C}, \utilde{A}_\mu, \utilde{A}_{\mu, \alpha} (\alpha < \mu)$
are also $\mathbb{P}_\vartheta$-names. By Proposition \ref{dicht} (2) we have
$\Vdash_{\mathbb{P}_\vartheta}\, ''\utilde{A}_\mu$ is an a.d. family$''$. In order to prove
maximality, and hence $\Vdash_{\mathbb{P}_\vartheta} \mu \in \utilde{\mathcal{A}}$, we
proceed completely analogously to Proposition \ref{dicht} (2A):

\begin{cla}\label{ka}
${\Vdash_{\mathbb{P}_\vartheta}}\, ''\utilde{A}_\mu$ is a mad familiy $''$.
\end{cla}

\textbf{Proof of Claim \ref{ka}:} By contradiction assume that $p \in \mathbb{P}_\vartheta$
and $\utilde{a}$ is a $\mathbb{P}_\vartheta$-name such that ${p \Vdash_{\mathbb{P}_\vartheta}}\,
'' \utilde{a} \in [\omega]^\omega$ and $\utilde{a} \notin \utilde{A}_\mu$ and
$\utilde{A}_\mu \cup \{\utilde{a}\}$ is a.d.$''$. By Claim \ref{fix} there is
$f \in J_\vartheta$ such that $p \in \mathbb{P}_{\vartheta, f}$ and $\utilde{a}$ is a
$\mathbb{P}_{\vartheta, f}$-name. W. l. o. g. we may assume that $\mu \in {\rm dom}\,
(f (0))$. As $f (0) (\mu) \in [\mu]^{<\vartheta}$ and $\vartheta \leqslant \mu$, we can
choose $\alpha \in \mu \setminus f (0) (\mu)$. We can find $q \geqslant_{\mathbb{P}_\vartheta}
p$ and $m < \omega$ such that $q \Vdash_{\mathbb{P}_\vartheta} \utilde{A}_{\mu, \alpha}
\cap \utilde{a} \subseteq m$. Choose $p_1 \in \mathbb{P}_{\vartheta, f}, p_1 \geqslant q
\upharpoonright f$ (see (K)), and $k\geqslant n_{q (0) (\mu)}$
(see \ref{hei} (A)) such that $p_1 {\Vdash_{\mathbb{P}_{\vartheta, f}}}\, '' k \in
\utilde{a} \setminus \utilde{a}_{\mu, \beta}\,''$ for all $\beta \in f (0) (m) \cap
{\rm dom}\, (q (0) (\mu))$. Similarly to \ref{dicht} (2A) we can define $q_1 \in
\mathbb{P}_\vartheta$ such that $q_1 \geqslant q, q_1 \upharpoonright f \geqslant p_1$ and
$q_1 \Vdash_{\mathbb{P}_\vartheta} k \in \utilde{a} \cap \utilde{a}_{\mu, \alpha}$,
which is a contradiction. \hfill
$\Box_{\ref{ka}}$

To prove $\mathop{\Vdash}_{\mathbb{P}_\mathbb{C}} \min (\utilde{\mathcal{A}}) = \min
(\mathbb{C})$ we have to recall that $\mathbb{P}_\vartheta / \mathbb{Q}_0$ forces
$MA_{<\vartheta}$ (Knaster). Moreover, given an a. d. family $A = \langle a_\alpha :
\alpha < \mu\rangle, \mu \geqslant \omega$, there exists a standard $\sigma$-centered
forcing notion $Q_A$ which adds $a \in [\omega]^\omega$ such that
$A \cup \{a\}$ is a. d. Its conditions are pairs $(x, F) \in [\omega]^{<\omega} \times
[A]^{<\omega}$ ordered by $(x, F) \leqslant (y, H)$ iff $x \subseteq y, F \subseteq H$
and $y \setminus x \cap a_\alpha = \emptyset$ for every $a_\alpha \in F$. Now if
$\mathbb{P}_\vartheta$ added some mad familiy $A$ of size $\omega \leqslant \mu < \vartheta$,
it hat to be added by $\mathbb{P}_\alpha$ for some $\alpha < \vartheta$. But then one
of the factors of $\utilde{\mathbb{Q}}_\alpha$ is an isomorphic copy of $Q_A$ (see
(H), (I)), and hence $A$ is not maximal after forcing with $\mathbb{P}_{\alpha +1}$.

It remains to prove that after forcing with $\mathbb{P}_\vartheta$ no cardinal $\lambda
\in [\min (\mathbb{C}), \max (\mathbb{C})] \setminus \mathbb{C}$ belongs to the mad
spectrum. For this we shall generalize the arguments from \S\, 2. Let $\sigma$ be the
minimal regular cardinal $\leqslant \lambda$ such that $[\sigma, \lambda] \cap
\mathbb{C} = \emptyset$.

Towards a contradiction assume ${p \mathop{\Vdash}_{\mathbb{P}_\vartheta}} ''\langle
\utilde{B}_\alpha : \alpha < \lambda\rangle$ is a m.a.d. family$''$.

By Claim \ref{fix} for each $\alpha < \lambda$ we have the following:

\begin{itemize}
\item[$(\ast)_7$]
\begin{itemize}
\item[(a)] $f_\alpha \in J_\vartheta$ such that $0 \in {\rm dom}\, (f_\alpha),
p \in \mathbb{P}_{\vartheta, f_\alpha}$ and $\utilde{B}_\alpha$ is a $\mathbb{P}_{\vartheta,
f_\alpha}$-name;
\item[(b)] $\mathbb{B}_\alpha$ is a Borel function such that $\utilde{B}_\alpha = \mathbb{B}_\alpha
(\ldots, \utilde{\eta}_{\mu, (\alpha, n), \xi (\alpha, n, m)}, \ldots,\linebreak \utilde{\nu}_{j (\alpha,
n), \varepsilon (\alpha, n, m)} (\gamma (\alpha, n, m)) \ldots)_{ n, m}$, where
\item[(c)] $\langle \mu (\alpha, n) : n < \omega\rangle$ is with no repetition, $\mu (\alpha, n) \in
{\rm dom}\, (f_\alpha (0)) (\subseteq \mathbb{C})$;
\item[(d)] $\langle \xi (\alpha, n, m) : m < \omega\rangle$ is with no repetition, $\xi (\alpha, n, m)
\in f_\alpha (0) (\mu (\alpha, n))$;
\item[(e)] $\langle j (\alpha, n) : n < \omega\rangle$ is with no repetition, $j (\alpha, n) \in
{\rm dom}\, (f_\alpha) \setminus \{0\}\, \, (\subseteq \vartheta)$;
\item[(f)] $\langle \varepsilon (\alpha, n, m) : m < \omega\rangle$ is with no repetition,
$\varepsilon (\alpha, n, m) \in f_\alpha (j (\alpha, n))$;
\item[(g)] $\gamma (\alpha, n, m) < \gamma_{j (\alpha, n), \varepsilon (\alpha, n, m)}$ (which is
the domain of $\utilde{\mathbb{Q}}_{j (\alpha, n), \varepsilon (\alpha, n, m)})$.
\end{itemize}
\end{itemize}

\begin{defn}
We define a binary relation $E$ on $\sigma$ by letting $(\alpha, \beta) \in E$ iff

\begin{itemize}
\item[(a)] there exists $\bm{s}_{\alpha, \beta} \in {\rm paut}_\alpha$ such that $\bm{s}_{\alpha,
\beta} = (f_\alpha, f_\beta, \mathbf{x}_{\alpha, \beta}, \overline{k}^{\alpha, \beta})$, hence
${\rm dom}\, (f_\alpha) = {\rm dom}\, (f_\beta)$ and $\kappa_{\alpha, \beta} :=
\kappa_{\bm{s}_{\alpha, \beta}}$ is an isomorphism from $\mathbb{P}_{\vartheta, f_\alpha}$ onto
$\mathbb{P}_{\vartheta, f_\beta}$ in particular, and $k_i^{\alpha, \beta} : f_\alpha (i) \rightarrow
f_\beta (i)$ for $i \in {\rm dom}\, f_\alpha \setminus \{0\}$ is order-preserving, hence o. t.
$f_\alpha (i) = o.t.\, f_\beta (i)$, and
\item[(b)] the isomorphism $\kappa_{\alpha, \beta}$ from (a) maps the $\mathbb{P}_{\vartheta, f_\alpha}$-name
$\utilde{B}_\alpha$ onto the $\mathbb{P}_{\vartheta, f_\beta}$-name $\utilde{B}_\beta$.
\end{itemize}
\end{defn}

It is clear that $E$ is an equivalence relation. Note that $E$ has no more than $\vartheta$ many equivalence
classes. Indeed, for given $\alpha < \sigma$ we can recursively define $g \in J_{i^*}$, where
$i^* = \sup ({\rm dom}\, (f_\alpha))$, with ${\rm dom}\, (g) = {\rm dom}\, (f_\alpha)$, some finite
support iteration $\langle \mathbb{P}_{i, g}, \utilde{\mathbb{Q}}_j^g : i \in {\rm dom}\, (g) \cup
\{i^*\}, j \in {\rm dom}\, (g) \rangle$ and $\langle \kappa_i : i \in {\rm dom}\, (g) \cup \{i^*\}\rangle$ such
that

\begin{itemize}
\item[(a)] $g (0) = \langle u_\mu : \mu \in {\rm dom}\, (g (0))\rangle$ for some ${\rm dom}\, (g (0))
\subseteq \vartheta$ and $u_\mu \subseteq \vartheta$ with $\mathop{\Sigma}\limits_{\mu \in {\rm dom}\, (g (0))}
|u_\mu| + 1 < \vartheta$ such that, letting $f_\alpha (0) = \langle v_\mu : \mu \in {\rm dom}\, (f_\alpha)\rangle$,
we have $|{\rm dom}\, g (0)| = |{\rm dom}\, (f_\alpha (0))|$ and some bijection $\pi : {\rm dom}\, (g (0))
\rightarrow {\rm dom}\, (f_\alpha (0))$ such that $|u_\mu| = |v_{\pi (\mu)}|$ for every $\mu \in {\rm dom}\, (g)$,
\item[(b)] $g (i) = o.t.\, (f_\alpha (i))$ for $i \in {\rm dom}\, (g) \setminus \{0\}$,
\item[(c)] for every $j \in {\rm dom}\, (g)\, \, {\mathop{\Vdash}_{\mathbb{P}_{j, g}}} '' \utilde{\mathbb{Q}}_j^g =
\prod\limits_{\nu < g (i)} \utilde{\mathbb{Q}}^g_{j, \nu}$ is a finite support product where $\utilde{\mathbb{Q}}^g_{j,
\nu}$ has the Knaster property and ${\rm dom}\, (\utilde{\mathbb{Q}}^g_{j, \nu}) = {\gamma_{j, \varepsilon_\nu}}''$,
where $\varepsilon_\nu$ is the $\nu$-th element of $f_\alpha (j)$, and $\kappa_j : \mathbb{P}_{j, g} \rightarrow
\mathbb{P}_{j, f_\alpha}$ is an isomorphism, such that for every $p \in \mathbb{P}_{j, g}$ and $\xi, \zeta \in
\gamma_{i, \nu}$ we have that $p\, {\mathop{\Vdash}_{\mathbb{P}_{j, g}}} '' \xi \mathop{<}_{\utilde{\mathbb{Q}}^g_{j, \nu}}
{\zeta} ''$ iff $\kappa_j (p) {\mathop{\Vdash}_{\mathbb{P}_{j, f_\alpha}}} ''\xi \mathop{<}_{\utilde{\mathbb{Q}}_{j, \nu}}
{\zeta} ''$. Hence $\kappa_j$ can be extended in a natural way to an isomorphism between $\mathbb{P}_{j, g} *
\utilde{\mathbb{Q}}^g_j$ and $\mathbb{P}_{j, f_\alpha} * \utilde{\mathbb{Q}}_j$.

\end{itemize}

Finally let $\utilde{C}$ be the $\mathbb{P}_{i^*, g}$-name that by $\kappa_{i^*}$ is mapped to $\utilde{B}_\alpha$.

By our assumption $\vartheta = \vartheta^{<\vartheta}$ and by property $(J)$ it is clear that there
are at most $\vartheta$ many $g, \langle\mathbb{P}_{i, g}, \utilde{\mathbb{Q}}_j^g : i \in {\rm dom}\, (g) \cup
\{i^*\}, j \in {\rm dom}\, (g)\rangle$ and $\utilde{C}$ as above. Moreover, if $\alpha, \beta < \sigma$
produce the same these objects, then $\alpha\, E\, \beta$.

Therefore, without loss of generality we may assume $\alpha\, E\, \beta$ for all $\alpha, \beta < \sigma$. In
particular $\mathbb{B}_\alpha = \mathbb{B}_*$ and $j (\alpha, n) = j (n)$ for all $\alpha < \sigma, n < \omega$.
Similarly to the proof of \ref{golo} we shall now normalize the remaining relevant indices in the computation
$(*)_7$ (b) of $\utilde{B}_\alpha$. Actually it is more convenient to normalize the $f_\alpha$.

As we assume $\alpha\, E\, \beta$ for every $\alpha, \beta < \sigma$, if $f_\alpha (0) = \langle u^\alpha_\nu
: \nu \in C^\alpha\rangle$, then $\delta (0) := |C^\alpha |$ and $\delta (1 + \nu) := |u^\alpha_\nu|$ for
$\nu < \delta (0)$ do not depend on $\alpha$. Let $\langle \mu (\alpha, \nu) : \nu < \delta (0)\rangle$
enumerate $C^\alpha$ and $\langle \xi (\alpha, \nu, \mu) : \mu < \delta  (1 + \nu)\rangle$ enumerate $u^\alpha_\nu$.
Similarly, $d^* := {\rm dom}\, (f_\alpha)$ and $o_\nu := o.t.\, f_\alpha (\nu)$ do not depend on
$\alpha < \sigma$ for $\nu \in d^* \setminus \{0\}$. Let $\langle\varepsilon (\alpha, \nu, \mu) : \mu <
o_\nu\rangle$ increasingly enumerate $f_\alpha (\nu)$. We let $\delta (\delta (0) + \mu) = o_{\nu_{1+\mu}}$
for every $\mu < {\rm o.\, t.}\, (d^*) =: o^*$, where $\nu_\mu$ is the $\mu$-th element of $d^*$.

Now for each $\alpha < \sigma$ we define $\overline{\zeta}_\alpha = \langle\zeta_\alpha (i) : i <
\mathop{\Sigma}\limits_{\nu < \delta (0) + o^*} \delta (\nu)\rangle$ (Without loss of generality we assume that
$o^*$ is a limit ordinal.), such that

\begin{itemize}
\item[(a)] $\zeta_\alpha (i) = \mu (\alpha, i)$ for $i < \delta (0)$,
\item[(b)] $\zeta_\alpha ((\mathop{\Sigma}\limits_{\mu < 1 + \nu} \delta (\mu)) + i) = \xi
(\alpha, \nu, i)$ for $\nu < \delta (0)$ and $i < \delta (1 + \nu)$,
\item[(c)] $\zeta_\alpha ((\mathop{\Sigma}\limits_{\mu < \delta (0) + \nu} \delta (\mu)) + i)
= \varepsilon (\alpha, \nu, i)$ for $\nu < o^*$ and $i < \delta (\delta (0) + \nu)$.
\end{itemize}

Analogously to \ref{golo} we distinguish Cases A and B, where now $\chi = \min \{\mu : \mu^{<\vartheta}
\geqslant \sigma\}$ and

\textbf{Case A:} $\chi = \sigma$,\vspace{0.5ex}\\
\textbf{Case B:} $\chi < \sigma$, hence $\chi > \vartheta$ and $cf (\chi) < \vartheta$.

As there, applying the usual pigeon-hole principle in Case A, and [Sh620, 7.1 (0), (1)] in Case B
with $\lambda, \kappa, \mu, D, \langle f_\alpha : \alpha < \lambda\rangle$ there standing for
$\sigma, \mathop{\Sigma}\limits_{\nu < \delta (0) + o^*}, \vartheta^+, D_\sigma^{cb},
\langle \overline{\zeta}_\alpha : \alpha < \sigma\rangle$ here, without loss of generality we
may assume that for some partition $\mathop{\Sigma}\limits_{\nu < \delta (0) + o^*} \delta (\nu)
= w_0\, \dot\cup\, w_1\, \dot\cup\, w_2$ ($w_1 = \emptyset$ in Case A) and ordinal function $\overline{\beta}^*
=\langle\beta^* (i) : i \in w_0 \cup w_1\rangle$ we have

\begin{itemize}
\item[$(\ast)_8$]
\begin{itemize}
\item[(a)] cf $(\beta^* (i)) > \vartheta$ for every $i \in w_1$;
\item[(b)] for every $\alpha < \sigma$
\begin{itemize}
\item[($\alpha$)] $\overline{\zeta}_\alpha \upharpoonright w_0 = \overline{\beta}^*
\upharpoonright w_0$,
\item[($\beta$)] $\overline{\zeta}_\alpha \upharpoonright w_1 < \overline{\beta}^*
\upharpoonright w_1$;
\end{itemize}
\item[(c)] if $\overline{\gamma} \in \prod \overline{\beta}^* \upharpoonright w_1$, then for
$\sigma$ many $\alpha < \sigma$ we have $\overline{\gamma} < \overline{\zeta}_\alpha
\upharpoonright w_1$;
\item[(d)] for every $\alpha < \sigma$ and $i \in w_2\, \, \, \zeta_{\alpha} (i) \notin \{\zeta_\beta
(j) : \beta < \alpha, j < \mathop{\Sigma}\limits_{\nu < \delta (0) + o^*} \delta (\nu)\}$.
\end{itemize}
\end{itemize}

Note that $\beta^* \upharpoonright w_0$ is essentially some $f_* \in J_\vartheta$
with $f^* \leqslant f_\alpha$ for all $\alpha < \sigma$. Indeed let ${\rm dom}\, f_*$ contain
$0$ iff $w_0 \cap \delta (0) \neq \emptyset$, and contain $\nu \in d^* \setminus \{0\}$ iff
$w_0 \cap [\mathop{\Sigma}\limits_{\mu < \delta (0) + \nu} \delta (\mu),
\mathop{\Sigma}\limits_{\mu \leqslant \delta (0) + \nu} \delta (\mu)) \neq \emptyset$. Note
that if $\varepsilon$ belongs to this last intersection, then $g_{\nu, \varepsilon}
\leqslant_{J_\vartheta} f_\alpha \upharpoonright \nu$, hence $g_{\nu, \varepsilon} (0)
\leqslant_{I_0} f_\alpha (0)$ for all $\alpha$ and therefore $0 \in {\rm dom}\, f_*$.
If ${\rm dom}\, f_* \neq \emptyset$ and hence $0 \in {\rm dom}\, f_*$, we let ${\rm dom}\, f_* (0)
= \{\beta^* (i) : i \in \delta (0) \cap w_0\}$ and for $\beta^* (i) = \mu \in {\rm dom}\, f_* (0)$
we let $f_* (0) (\mu) = \{\beta^* (\mathop{\Sigma}\limits_{\nu \leqslant \delta (0) + i}
\delta (\nu) + j) : j < \delta (1 + i), \mathop{\Sigma}\limits_{\nu < \delta (0) + i}
+ j \in w_0\}$. For $\nu \in {\rm dom}\, f_* \setminus\{0\}$ we let $f_* (\nu) =
\{\beta^* ((\mathop{\Sigma}\limits_{\mu < \delta (0) + \nu} \delta (\mu) + i) : i < \delta
(\delta (0) + \nu), (\mathop{\Sigma}\limits_{\mu < \delta (0) + \nu} \delta (\mu)) + i \in w_0\}$.
The argument with $g_{\nu, \varepsilon}$ given above shows that $f_* \in J_\vartheta$.
Clearly $f_* \leqslant f_\alpha$, for all $\alpha < \sigma$.

Recursively we shall define $f_\lambda \in J_\vartheta$ and $\bm{s}^{\nu, \alpha}
\in {\rm paut}_\nu$ for $\nu \in {\rm dom}\, (f_\lambda) \cup \{\sup {\rm dom}\, (f_\lambda)\}$
and $\alpha < \sigma$ such that

\begin{itemize}
\item[$(\ast)_9$]
\begin{itemize}
\item[(a)] ${\rm dom}\, (f_\lambda) = d_*$,
\begin{itemize}
\item[($\alpha$)] $\langle \mu (\lambda, \nu) : \nu < \delta (0)\rangle$ enumerates
${\rm dom}\, (f_\lambda (0))$ and $\langle \xi (\lambda, \nu, \mu) : \mu < \delta
(1+\nu)\rangle$ enumerates $u^\lambda_\nu := f_\lambda (0) (\nu)$ for $\nu < \delta (0)$,\vspace{0.5ex}
\item[($\beta$)] $\langle\varepsilon (\lambda, \nu, i) : i < o_\nu\rangle$ enumerates
$f_\lambda (\nu)$;
\end{itemize}
\item[(b)] letting $\kappa^\nu_\alpha = \kappa_{\bm{s}^{\nu, \alpha}}$, we have that $\kappa^\nu_\alpha$
is an isomorphism from $\mathbb{P}_{\vartheta, f_\alpha}$ onto $\mathbb{P}_{\vartheta, f_\lambda}$
such that $\kappa_\alpha^\nu = \kappa^\nu_\beta \circ \kappa_{\alpha, \beta}$ for all
$\alpha < \beta < \sigma$ (here $\kappa_{\alpha, \beta}$ is the isomorphism from
$\mathbb{P}_{\vartheta, f_\alpha}$ onto $\mathbb{P}_{\vartheta, f_\beta}$ witnessing
$\alpha\, E\, \beta$).
\end{itemize}
\end{itemize}

We define $f_\lambda (0)$ analogously to $(*)_3$:

\begin{itemize}
\item[$(\ast)_{10}$]
\begin{itemize}
\item[(a)] If $\nu \in \delta (0) \cap w_0$, then $\mu (\lambda, \nu) = \beta^* (\nu)$;
\item[(b)] if $\nu \in \delta (0) \setminus w_0$, then $\mu (\lambda, n) = \Upsilon_{i (\lambda, \nu)}$,
where $i (\lambda, \nu)$ is the $\nu$-th member of $\{i < \max (\mathbb{C}) : \Upsilon_i \notin
\{\mu (\alpha, \nu) : \alpha < \lambda, \nu < \delta (0)\}\}$;
\item[(c)] if $\mathop{\Sigma}\limits_{\rho < 1 + \nu} \delta (\rho) + \mu \in w_0$, then
$\xi (\lambda, \nu, \mu) = \beta^* (\mathop{\Sigma}\limits_{\rho < 1 + \nu} \delta (\rho) + \mu)$;
\item[(d)] if $\nu \in \delta (0) \cap w_0$ and $\mathop{\Sigma}\limits_{\rho < 1 + \nu} \delta (\rho)
+ \mu \notin w_0$, then $\xi (\lambda, \nu, \mu)$ is the $\mu$-th member of $\mu (\lambda, \nu)
\setminus \{\xi (\beta, \nu_1, \mu_1) : \beta < \lambda, \nu_1 < \delta (0), \mu_1 < \delta (1 + \nu)\}$
(as in $(*)_3$ (d) this choice is possible);
\item[(e)] if $\nu \notin \delta (0) \cap w_0$ and $\mathop{\Sigma}\limits_{\rho < 1 + \nu} \delta (\rho)
+ \mu \notin w_0$, then $\xi (\lambda, \nu, \mu) = \mu$.
\end{itemize}
\end{itemize}

For each $\alpha < \sigma$ we have $\mathbf{x}_\alpha \in {\rm paut}_\mathbb{C}$ such that
$\mathbf{x}_\alpha = \langle g_\alpha, \overline{h}^\alpha, {\rm dom}\, (f_\alpha (0)), f_\alpha (0),
\linebreak {\rm dom}\, (f_\lambda (0)), f_\lambda (0)\rangle, g_\alpha (\mu (\alpha, \nu)) = \mu (\lambda, \nu)$
for $\nu < \delta (0), \overline{h}^\alpha = \langle h^\alpha_\nu :\nu < \delta (0)\rangle$ and $h^\alpha_\nu
(\xi (\alpha, \nu, \mu)) = \xi (\lambda, \nu, \mu)$ for $\mu < \delta (1 + \nu)$. Then $\kappa^1_\alpha
:= \kappa_{\mathbf{x}_\alpha}$ is an isomorphism between $\mathbb{Q}_{f_\alpha (0)}$ and
$\mathbb{Q}_{f_\lambda (0)}$.

Similarly we have isomorphisms $\kappa_{\alpha \beta}$ between $\mathbb{Q}_{f_\alpha (0)}$ and
$\mathbb{Q}_{f_\beta (0)}$ for $\alpha < \beta < \sigma$ such that $\kappa^1_\alpha = \kappa^1_\beta \circ
\kappa_{\alpha \beta}$.

Now suppose that $\nu \in d^*$ and we have constructed $\bm{s}^{\nu, \alpha} \in {\rm paut}_\nu$
with induced isomorphism $\kappa^\nu_\alpha  := \kappa_{\bm{s}^{\nu, \alpha}}$ from $\mathbb{P}_{\nu, f_\alpha
\upharpoonright \nu}$ onto $\mathbb{P}_{\nu, f_\lambda \upharpoonright \nu}$ for every $\alpha <
\sigma$ such that $\kappa^\nu_\alpha = \kappa_\beta^\nu \circ \kappa_{\alpha, \beta}$ for all
$\alpha < \beta < \sigma$. Suppose $\nu = \nu_\mu$ is the $\mu$-th element of $d^*$. We shall define
$\langle\varepsilon (\lambda, \nu, i) : i < o_\nu\rangle$ and let $f_\lambda (\nu) = \{\varepsilon
(\lambda, \nu, i) : i < o_\nu\}$. Then we define $k^\alpha_\nu (\varepsilon (\alpha, \nu, i)) = \varepsilon
(\lambda, \nu, i)$ for every $\alpha < \sigma$ and $i < o_\nu$, so that $k^\alpha_\nu$ extends $\bm{s}^{\nu, \alpha}$
as desired.

Suppose that $\varepsilon (\lambda, \nu, j)$ have been defined for all $j < i$. We have to define
$\varepsilon := \varepsilon (\lambda, \nu, i)$ in such a way that the two demands in Definition \ref{logo}
are satisfied: Firstly, $g_{\nu, \varepsilon (\alpha, \nu, i)}$ is mapped to $g_{\nu, \varepsilon}$
by $\kappa^\nu_\alpha$ and, secondly, the pair $\kappa^\nu_\alpha, id_{\gamma_{\nu, \varepsilon (\alpha,
\nu, i)}}$ maps $\utilde{\mathbb{Q}}_{\nu, \varepsilon (\alpha, \nu, i)}$ onto $\utilde{\mathbb{Q}}_{\nu,
\varepsilon}$, for every $\alpha < \sigma$. In case $\mathop{\Sigma}\limits_{\rho < \delta (0) + \mu}
\delta (\rho) + i \in w_0$, we let $\varepsilon (\lambda, \nu, i) = \beta^* (\mathop{\Sigma}\limits_{\rho
< \delta (0) + \mu} \delta (\rho) + i)$. Then clearly these demands are satisfied.

Now suppose $\mathop{\Sigma}\limits_{\rho < \delta (0) + \mu} \delta (\rho) + i \notin w_0$. By
construction we have that $k^{\alpha, \beta}_\nu (\varepsilon (\alpha, \nu, i))\linebreak = \varepsilon
(\beta, \nu, i), \kappa_{\alpha, \beta}$ maps $g_{\nu, \varepsilon (\alpha, \nu, i)}$ to
$g_{\nu, \varepsilon (\beta, \nu, i)}, \gamma := \gamma_{\nu, \varepsilon (\alpha, \nu, i)} =
\gamma_{\nu, \varepsilon (\beta, \nu, i)}$ and the pair $\kappa_{\alpha, \beta}, id_\gamma$ maps
$\utilde{\mathbb{Q}}_{\nu, \varepsilon (\alpha, \nu, i)}$ to $\utilde{\mathbb{Q}}_{\nu, \varepsilon
(\beta, \nu, i)}$ for all $\alpha < \beta < \sigma$.

Since $\kappa^\nu_\alpha$ and $\kappa_{\alpha, \beta}$ commute, we have that, letting $g:= \kappa^\nu_\alpha
(g_{\nu, \varepsilon (\alpha, \nu, i)})$ and $\utilde{\mathbb{Q}}$ the image of $\utilde{\mathbb{Q}}_{\nu,
\varepsilon (\alpha, \nu, i)}, g$ and $\utilde{\mathbb{Q}}$ do not depend on $\alpha$. Applying (H) ($\gamma$)
we choose $\varepsilon$ minimal in $\max (\mathbb{C}) \setminus (\bigcup \{f_\alpha (\nu) : \alpha
< \lambda\} \cup \{\varepsilon (\lambda, \nu, j) : j < i\})$ such that $\langle \utilde{\mathbb{Q}}_{\nu, \varepsilon},
g_{\nu, \varepsilon}\rangle = \langle \mathbb{Q}, g\rangle$. Defining $\varepsilon (\lambda, \nu, i) := \varepsilon$,
we can easily verify that all demands are satisfied.

If $\nu \in d^*$ is such that $d^* \cap \nu$ has no maximum or $\nu = \sup d^*$, we define $\kappa^\nu_\alpha =
\bigcup\limits_{\mu \in d^* \cap \nu} \kappa^\mu_\alpha$. Then $(*)_9$ holds. Let $\nu^* = \sup d^*$.
We define $\overline{\zeta}_\lambda$ for $f_\lambda$ precisely as we defined $\overline{\zeta}_\alpha$ for $f_\alpha,
\alpha < \sigma$ above.

Fix $\alpha < \sigma$ and define $\utilde{B}_\lambda = \kappa^{\nu^*}_\alpha (\utilde{B}_\alpha)$. By
the definition of the relation $E$ and by $(*)_9$ (b) we have that $\utilde{B}_\lambda$ does not depend on
$\alpha$.

Analogously to $(*)_1$ we can now show $p \Vdash_{\mathbb{P}_\vartheta}\, ''\langle \utilde{B}_\alpha :
\alpha \leqslant \lambda \rangle$ is an a. d. family$''$ and thus reach our desired contradiction. Indeed, let
$\alpha < \lambda$ be arbitrary. By $(*)_8$ we can find $\beta \in \sigma \setminus \{\alpha\}$ so that outside
$f_*, f_\beta$ and $f_\alpha$ are disjoint, i. e. if $i \in \mathop{\Sigma}\limits_{\nu < \delta (0) + o^*}
\delta (\nu) \setminus w_0$ then $\zeta_\beta (i) \notin \bigcup \{f_\alpha (\nu) : \nu \in {\rm dom}\, (f_\alpha)
\} \cup \bigcup \{f_\alpha (0) (\mu) : \mu \in {\rm dom}\, f_\alpha (0)\}$.

Hence we can define a bijection $\pi$ between

\[\begin{array}{l}\bigcup\, \{f_\alpha (\nu) \cup f_\beta (\nu) : \nu \in {\rm dom}\, (f_\alpha)
\cup {\rm dom}\, (f_\beta)\}\, \cup\vspace{2ex} \\ \bigcup\, \{f_\alpha (0) (\mu) \cup f_\beta (0) (\mu) : \mu \in {\rm dom}\, (f_\alpha (0))
\cup {\rm dom}\, (f_\beta (0))\}\end{array}\]

and

\[\begin{array}{l}\bigcup\, \{f_\alpha (\nu) \cup f_\lambda (\nu) : \nu \in {\rm dom}\, (f_\alpha)
\cup {\rm dom}\, (f_\lambda)\}\, \cup\vspace{2ex}\\ \bigcup\, \{f_\alpha (0) (\mu) \cup f_\lambda (0) (\mu) : \mu \in {\rm dom}\, (f_\alpha
(0)) \cup {\rm dom}\, (f_\lambda (0))\}\end{array}\]

so that $\pi$ is the identity except for $\pi (\zeta_\beta (i)) = \zeta_\lambda (i)$
in case $i \in \mathop{\Sigma}\limits_{\nu < \delta (0) + o^*} \delta (\nu) \setminus w_0$. Then $\pi$
induces an isomorphism $\kappa$ between $\mathbb{P}_{\vartheta, f_\alpha \cup f_\beta}$ and $\mathbb{P}_{\vartheta,
f_\alpha \cup f_\lambda}$ which fixes $p$ and $\utilde{B}_\alpha$, but maps $\utilde{B}_\beta$
to $\utilde{B}_\lambda$. As in $(*)_1$ this is a contradiction.\hfill $\Box$

\vspace{3cm}

\textbf{References}

\begin{itemize}
\item[{[Br]}] J\"{o}rg Brendle, The almost-disjointness number may have countable cofinality, Transactions of the
American Mathematical Society 355 (2003), no. 7, 2633-2649 (electronic).

\item[{[H]}] Stephen S. Hechler, Short complete nested sequences in $\beta N \setminus N$ and small maximal
almost-disjoint families, General Topology and its Applications 2 (1972), 139-149.

\item[{[Sh620]}] Saharon Shelah, Special Subsets of $^{cf (\mu)} \mu$, Boolean Algebras and Maharam measure
Algebras, Topology and its Application 99 (1999), 135-235, 8th Prague Topological Symposium on General Topology
and its Relations to Modern Analysis and Algebra, Part II (1996). math. LO/9804156.

\item[{[Sh700]}] Saharon Shelah, Two cardinal invariants of the continuum $(\frak{d} < \frak{a})$ and
FS linearly ordered iterated forcing, Acta Mathematica 192 (2004), 187-233.
\end{itemize}
\vspace{3cm}

First author: Einstein Institute of Mathematics, Edmond J. Safra Campus, Givat Ram, The Hebrew University
of Jerusalem, Jerusalem, 91904, Israel,

and

Department of Mathematics, Hill Center - Busch Campus, Rutgers, The State University of New Jersey, 110 Frelinghuysen
Road, Piscataway, NJ 08854-8019, USA

Second author: Mathematisches Seminar der Christian-Albrechts-Universit\"at zu Kiel, Ludewig-Meyn-Stra{\ss}e 4,
24118 Kiel, Germany

\end{document}